\documentclass[12pt]{article}
\usepackage[utf8]{inputenc}
\usepackage{authblk}

\usepackage[numbers]{natbib}
\usepackage{amssymb}
\usepackage{array} 
\usepackage{siunitx}[=v2] 
\DeclareSIUnit{\gramdryweight}{\ensuremath{g_{CDW}}}
\usepackage[version=4]{mhchem} 
\usepackage{amsmath}  
\usepackage{bm}  
\usepackage{subcaption} 
\usepackage{mathtools} 
\usepackage{booktabs} 
\usepackage{verbatim} 
\usepackage{multirow} 
\usepackage{hyperref}

\usepackage{setspace}
\usepackage[margin=1in]{geometry}

\begin{document}
\title{Simultaneous design of fermentation and microbe}
\author[1]{Anita L. Ziegler}
\author[1]{Ashutosh Manchanda}
\author[1]{Marc-Daniel Stumm}
\author[2]{Lars M. Blank} 
\author[3,1,4,*]{Alexander Mitsos}

\affil[1]{Process Systems Engineering (AVT.SVT), RWTH Aachen University, 52074 Aachen, Germany}
\affil[2]{Institute of Applied Microbiology (iAMB), Aachen Biology and Biotechnology (ABBt), RWTH Aachen University, 52074 Aachen, Germany}
\affil[3]{JARA-ENERGY, 52056 Aachen, Germany}
\affil[4]{Institute of Energy and Climate Research, Energy Systems Engineering (IEK-10), Forschungszentrum Jülich GmbH, 52425 Jülich, Germany}
\affil[*]{Corresponding author: Alexander Mitsos, amitsos@alum.mit.edu}

\date{\today}

\maketitle

\section*{Abstract}

Constraint-based optimization of microbial strains and model-based bioprocess design have been used extensively to enhance yields in biotechnological processes. However, strain and process optimization are usually carried out in sequential steps, causing underperformance of the biotechnological process when scaling up to industrial fermentation conditions.
Herein, we propose the optimization formulation \emph{SimulKnock} that combines the optimization of a fermentation process with metabolic network design in a bilevel optimization program. The upper level maximizes space-time yield and includes mass balances of a continuous fermentation, while the lower level is based on flux balance analysis. SimulKnock predicts optimal gene deletions and finds the optimal trade-off between growth rate and product yield. Results of a case study with a genome-scale metabolic model of \textit{E. coli} indicate higher space-time yields than a sequential approach using OptKnock for almost all target products considered. By leveraging SimulKnock, we reduce the gap between strain and process optimization.

\section*{Keywords}
process optimization, computational strain design, constraint-based metabolic modeling, metabolic engineering

\section{Introduction}
Industrial microbiology promises product synthesis from renewable feedstocks representing a sustainable alternative to petro-chemical synthesis \citep{Clomburg.2017}. 
Often, these biotechnological products are produced by high-performing microbial strains that have been designed using metabolic engineering. 
Strain design is supported by computational methods to reduce the experimental effort \citep{Landon.2019}.
In the first place, constraint-based metabolic modeling formulations predict the influence of a given genetic modification on the microorganism \citep{Maia.2016}. 
The best-known formulation is flux balance analysis (FBA) \citep{Varma.1994, Orth.2010b}.
Using linear programming, FBA predicts the internal fluxes of a microorganism based on its genome-scale metabolic model (GEM). 
Variations of FBA refine the cellular objective of the organism \citep{Schuetz.2012} or better suit a genetically modified microorganism \citep{Segre.2002, Brochado.2012, Shlomi.2005}.
Constraint-based strain optimization formulations go one step further and predict targets for genetic modification \citep{ValderramaGomez.2017, Maia.2016}. 
The first formulation, OptKnock, was presented by \citet{Burgard.2003}.
OptKnock proposes optimal gene knockouts by solving a bilevel optimization problem. 
The upper level represents the bioengineering perspective to maximize the product yield on the substrate. 
The lower level, based on FBA, represents the microorganism with the cellular objective of maximizing biomass production. 
OptKnock has been extended, e.g., to account for worst-case predictions \citep{Tepper.2010}, insertion of genes \citep{Pharkya.2004}, and up and down regulation of genes \citep{Pharkya.2006}.
It was also modified to suit a genetically modified organism better \citep{Kim.2011}.
For strain design, the predicted modifications are experimentally tested and evaluated at a laboratory scale, and the best-performing strain is chosen. 

For industrial production, industrial-scale process conditions come into play. 
Process design involves the analysis of mass and energy balances an analyzing production costs, including up- and downstream processing \citep{Biegler.1997}.
In bioprocess design, a suitable kinetic represents the microorganism, e.g., its growth \citep{Chmiel.2018, Villadsen.2011}.
Computational process optimization supports process design by adjusting the process parameters to achieve maximal space-time yield or minimal production cost \citep{Villadsen.2011, Edgar.2001, Gordeeva.2015}. 

Typically, strain and process design are separate steps, performed in sequence. 
Process conditions, however, influence the behavior and performance of the microorganism.
Hence, a high-performing strain in the laboratory may underperform in an industrial process. 
For example, a strain is first designed under small-scale batch fermentation conditions in the laboratory. 
Next, the industrial process is designed for the designated strain, whereby industrial processes are envisaged to run on a large scale in continuous mode, and example processes already exist in the pharmaceutical and food sectors.
This switch in conditions entails a significant adjustment for the organism.
Reasons for the under-performance can be steep gradients in the spatial distribution of substrates (e.g., oxygen) in the fermenter, high downstream processing costs, or lacking genetic and phenotypic stability by the microbe \citep{Wehrs.2019, Olsson.2022}. 
Strain and process design coupling is needed to overcome these scale-up difficulties \citep{Olsson.2022,Richelle.2020}.
There exist different approaches to do so. 
In their attempt to capture spatial changes in a stirred bioreactor, \citet{Lapin.2004} coupled computational fluid dynamics and metabolic modeling. 
To model diauxic growth, \citet{Mahadevan.2002} connected the mass balance equations of batch fermentation with metabolic modeling in their dynamic flux balance analysis (dFBA). 
Structurally, dFBA are differential‐algebraic equations with embedded optimization criteria.
For fed-batch fermentation, \citet{Oliveira.2021} replaced the metabolic network in dFBA with a surrogate model.
In an iterative process, \citet{Zhuang.2013} connected dFBA with computational strain optimization techniques to suit the organism to the future process conditions of batch fermentation. 
\citet{Ploch.2019} developed a differential-algebraic equation system with embedded dFBA optimization to model a biorefinery process under changing conditions. 
\citet{Jabarivelisdeh.2018} introduced a bilevel dynamic optimization framework, where a previously found genetic modification can be switched on during batch fermentation.
In a step-wise approach, \citet{TafurRangel.2022} selected a microbial strain regarding the required downstream process units. 
\citet{Dimitriou.2022} considered downstream process optimization simultaneous to strain optimization, namely in a superstructure optimization program.
However, the implications of the switch from laboratory (batch or
fed-batch) fermentation conditions to industrial (continuous) fermentation conditions have yet to be considered.

We propose a new formulation combining strain optimization with process optimization. 
Our simultaneous approach \emph{SimulKnock}, suggests optimal gene knockouts that suit a strain to its future industrial-scale conditions, and at the same time, the process conditions are customized to the strain. 
SimulKnock is a bilevel optimization program where the upper level maximizes the space-time yield subject to the mass balances of fermentation, and the lower level maximizes biomass production based on the FBA formulation.
In this paper, we consider continuous fermentation, an envisaged fermentation mode in industry. 
SimulKnock, however, can readily be extended to account for other fermentation modes, too. 
Certainly, extensions may result in different optimization classes, e.g., simultaneous optimization of batch operation and strain design would be a dynamic bilevel optimization problem, which is challenging \citep{Ploch.2019, Mitsos.2009}.
Michaelis-Menten or Monod kinetics are employed to connect the two levels.
Compared to existing strain optimization formulations, SimulKnock can be interpreted as an extension of OptKnock by process equations or as a variation of dFBA with a switch to continuous fermentation and implementation of a process objective and knockout predictions. 
Similarly, SimulKnock can be seen as an extension of process optimization by the strain optimization level (lower level).
For our case studies, we use strong dualization to reformulate the bilevel optimization program into a single-level program and solve it globally. 
We apply our formulation to a GEM of \textit{E. coli} for the production of ethanol, succinate, acetate, formate, fumarate, and lactate and compare the results with sequential optimization results and experimental results for succinate production from the literature. 
The proposed formulation SimulKnock will help close the gap in profitability between petrochemical and biochemical processes.

\section{Method: SimulKnock}
\subsection{Mass balances for a continuous process optimization} 

For optimal process design, the fermentation process was modeled with mathematical equations.
These equations may include mass and energy balances of different reactor types and process units. 
In SimulKnock, we chose the mass balance equations of a continuous stirred-tank reactor. 
In principle, however, SimulKnock can easily be extended to account for other reactor types or even full processes. 
We assumed that the reactor is ideally mixed and contains only one homogeneous phase, that only one substrate limits growth and production, and that only one population of organisms is in the tank. 
Furthermore, we assumed that the volumetric inflow equals the volumetric outflow. 
The biomass, the substrate, and the product were chosen as the representative compounds of the fermentation.  
Similarly to standard literature, mass conservation leads to
\begin{equation}
	\label{eqn:Conti_Process}
	\begin{aligned}
		\frac{\mathrm{d}c_{bio}}{\mathrm{d}t} &= (\mu-D)c_{bio}\\
		\frac{\mathrm{d}c_S}{\mathrm{d}t} &= -\frac{1}{Y_{bio/S}}\mu c_{bio} -\frac{1}{Y_{P/S}} q_p c_{bio} -m_S c_{bio}+D(c_{S,Feed}-c_S)\\
		\frac{\mathrm{d}c_{P}}{\mathrm{d}t} &= q_P c_{bio} -c_{P}\cdot D,
	\end{aligned}
\end{equation}
\\
where $D$ is the dilution rate in \si{\per\hour}; $c_{bio}$ is the biomass concentration in  gram Cell Dry Weight (CDW) per liter (\si{\gramdryweight\per\liter}); $c_{S,Feed}$, $c_S$ and $c_{P}$ are the substrate feed concentration, substrate concentration, and product concentration, respectively (all in \si{\gram\per\liter}). 
$Y_{bio/S}$ and $Y_{P/S}$, both with the unit \si{\gram\per\gram}, are the yields of biomass and product per substrate, respectively. 
The maintenance factor is denoted as $m_S$ in \si{\gram\per\gram\per\hour} and describes the required substrate uptake rate for cellular maintenance. 
Two kinetics were used in the equations: the product kinetics $q_P$ in \si{\gram\per\gram\per\hour}, describing the rate of formation of the product, and the growth rate kinetics $\mu$ in \si{\per\hour}, describing the rate of formation of the biomass.
The dilution rate $D$ is the quotient of the volumetric inflow/outflow and the culture volume within the bioreactor. 
The dilution rate is also the inverse of the residence time. 
As is typical, we assumed that a flow equilibrium is achieved after a sufficiently long period of operation at a constant dilution rate.
It follows that the continuous fermentation is at a steady state, such that $\frac{\mathrm{d}c_{bio}}{\mathrm{d}t} = 0$, $\frac{\mathrm{d}c_{S}}{\mathrm{d}t} = 0$, and $\frac{\mathrm{d}c_{P}}{\mathrm{d}t} = 0$ apply. 
Thereby, it directly follows that $D = \mu$, which means that biomass is neither washed out nor accumulated in the process.

\subsection{Flux balance analysis}
FBA analyzes the internal fluxes within a cell using linear programming \citep{Orth.2010b, Varma.1994}. 
The analysis is based on the metabolic network of an organism that includes the metabolites and the stoichiometry of the internal reactions. 
In FBA, based on the steady state assumption, the metabolites' mass balances are constrained to zero, and the reaction fluxes are constrained to upper and lower bounds.
One cellular objective is defined, e.g., the maximization of the biomass flux. 
The mathematical formulation of FBA reads
\begin{equation*}
	\begin{aligned}
		\max_{\bm{v} \in \mathbb{R}^n} \quad & v_{bio} \\
		\textrm{s.t.} \quad & \bm{Sv} = \bm{0} \\
		& \bm{v^{lower}} \le \bm{v} \le \bm{v^{upper}}
	\end{aligned}
\end{equation*}
where $\bm{v}$ denotes the vector of reaction fluxes in \si{\milli\mol\per\gramdryweight\per\hour} and where $n$ is the number of reactions in a metabolic network, where all reversible reactions were split into one forward and one backward reaction, further called irreversible network. 
The biomass flux is denoted by $v_{bio}$ in \si{\per\hour} and is an element of $\bm{v}$; $\bm{S}$ is the stoichiometric matrix of the irreversible network, and $lower$ and $upper$ are lower and upper bounds of the flux values, respectively. 
These bounds may include thresholds on the biomass flux or the ATP maintenance reaction. 

\subsection{Proposed combined optimization: SimulKnock} 

\begin{figure}[htb]
	\centering
	\includegraphics[trim = 3cm 5cm 6.5cmm 1cm, clip, width=\textwidth]{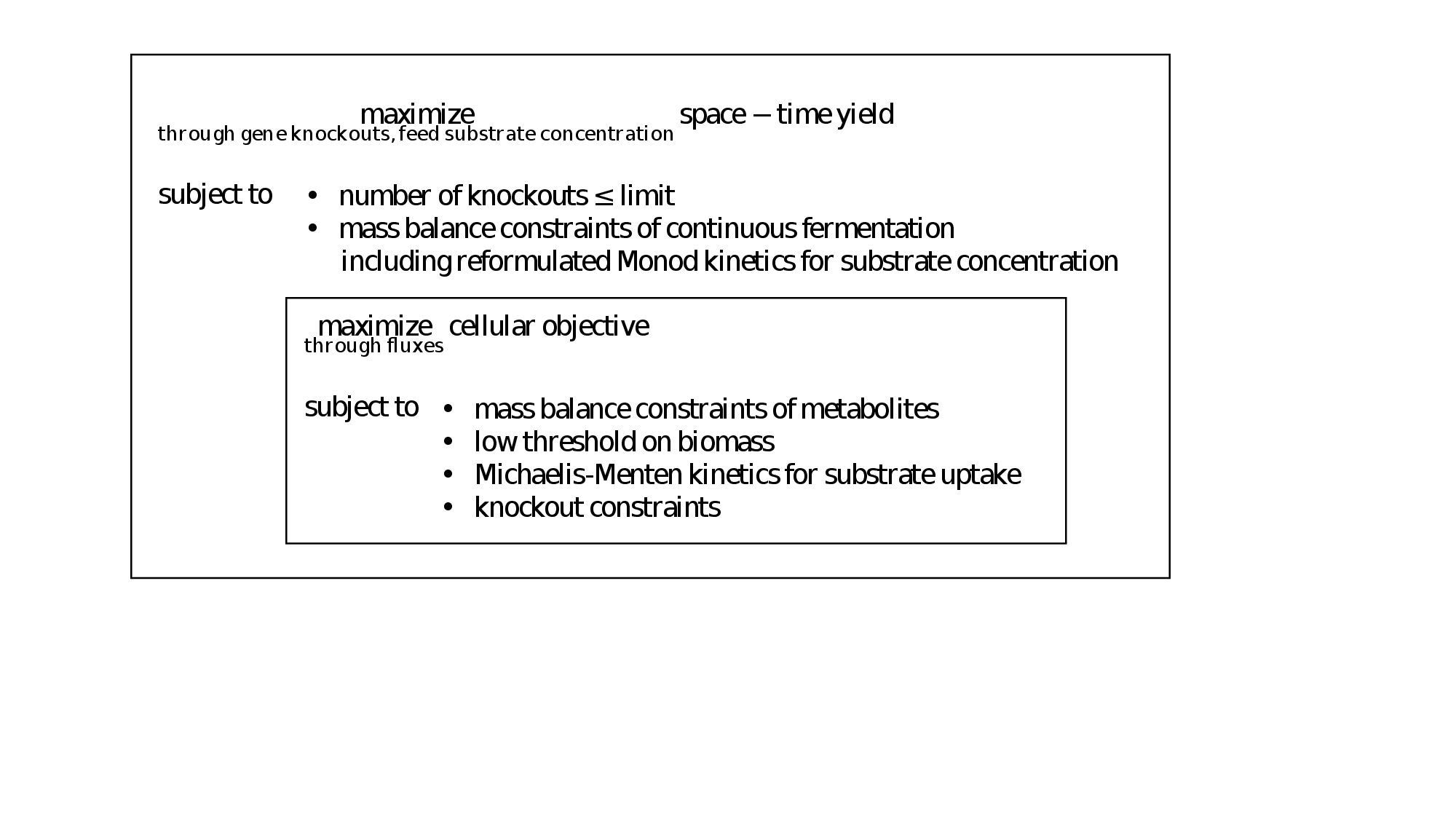}
	\caption{The bilevel optimization formulation of SimulKnock. Note that either Monod or Michaelis-Menten kinetics are applied at the upper level or the lower level. }
	\label{fig:Formulation_Boxes}
\end{figure}

We designed SimulKnock by embedding the FBA into the optimization of a continuous fermentation accounting for cellular degrees of freedom, resulting in a bilevel optimization formulation (Figure~\ref{fig:Formulation_Boxes}), following the paradigm of OptKnock \citep{Burgard.2003}. 
Using the steady-state mass balances of metabolites stemming from FBA, we assume that the cellular metabolism adapts infinitesimally quickly to shifting environmental conditions \citep{Stephanopoulos.1998}. 

SimulKnock maximizes the space-time yield in the upper-level program like other bioprocess optimization examples \citep{Gordeeva.2021, Sinner.2019}.
The connection between upper- and lower-level program is achieved via the expression of process parameters with metabolic variables, based on \citet{Ploch.2019} and \citet{Mahadevan.2002}. 
The growth rate directly transforms to $\mu = v_{bio}$, and due to steady-state continuous fermentation conditions, it holds $D=v_{bio}$.
Moreover, instead of considering that everything taken up goes into product formation, biomass formation, or maintenance, the substrate uptake is now directly expressed by $-\frac{1}{Y_{X/S}}\mu -\frac{1}{Y_{P/S}} q_p -m_S = -v_{S}\cdot M_S$, where $v_S$ is the substrate uptake flux in \si{\milli\mol\per\gramdryweight\per\liter} and $M_S$ denotes the molar mass of the substrate in \si{\gram\per\milli\mol}.
Thus, the maintenance is only directly considered at the lower level in SimulKnock by setting a threshold on the ATP maintenance reaction; there is no doubling at the upper level. 
The rate of product formation can also be directly considered now by setting $q_P = v_P \cdot M_P$, where $v_P$ is the product flux in \si{\milli\mol\per\gramdryweight\per\liter} and is an element of the flux vector $\bm{v}$ and $M_P$ is the molar mass of the product in \si{\gram\per\milli\mol}.
Hence, \eqref{eqn:Conti_Process} transforms to 
\begin{equation*}
	\begin{aligned}
		D&=v_{bio}\\
		0&=-v_{S}\cdot M_S \cdot c_{bio} +D\cdot(c_{S,f}-c_S)\\
		0&=v_{P}\cdot M_{P}\cdot c_{bio} -c_P\cdot D.
	\end{aligned}
\end{equation*}

Note that instead of fixing the dilution rate in the upper level and, thereby, fixing the growth rate, the dilution rate $D$ will be set according to the optimal value of $v_{bio}$ after the lower-level optimization. 
Hence, in the following, $D$ will be replaced with $v_{bio}$.

Furthermore, we considered two alternative kinetics to connect the upper and lower levels: Monod and Michaelis-Menten.
In Section~\ref{seq:inclusion_both_kinetics} in the Supporting Information, we discuss why including both kinetics simultaneously is not advisable.

Monod is an empirical, widely known, and easy-to-measure model of microbial growth. It links the substrate concentration with the growth rate and reads
\begin{align}
	v_{bio}&= v_{bio}^{max}\frac{c_S}{c_S+K_S} \xRightarrow[]{} c_S = K_S \frac{v_{bio}}{v_{bio}^{max}-v_{bio}}, &(Monod)
	\label{eq:Monod}
\end{align}
where $v_{bio}^{max}$ denotes the maximum growth rate in \si{\per\hour}, and $K_S$ denotes the Monod affinity constant  in \si{\gram\per\liter} for the specified substrate.
We assumed that gene deletions would not affect the kinetic parameters, i.e., we assumed that the kinetic parameters were constants for a given organism and substrate. 
Thus, both $v_{bio}^{max}$ and $K_S$ are parameters taken from literature or, in the case of $v_{bio}^{max}$, can be approximated with FBA. 
The growth rate results from the lower-level program, and thus, the Monod kinetics need to be imposed on the upper-level program (reactor level). 
As depicted on the right-hand side in \eqref{eq:Monod}, we, therefore, reformulated the kinetics such that the substrate concentration is a function of the growth rate: $c_S = f(v_{bio})$.  
Then, we implemented the kinetics in the upper level by replacing $c_S$ with the kinetic term. 
If we would keep $c_S$, the growth rate would be fixed by the kinetic term (equality constraint) instead of being optimized at the lower level. 
SimulKnock with included Monod kinetics reads

\begin{equation}
	\label{eq:SimulKnock_Monod}
	\begin{aligned}
		\underset{\boldsymbol{y},c_{S,Feed}, c_P, c_{bio}, v_S', v_{bio}', v_P'}\max &&\quad c_{P}\cdot v_{bio}'  &&          & \\
		\text{s.t. }                          &&\sum_{i=1}^{r} (1-y_i)&\leq K   & &\\
		&&0&=-v_{S}'\cdot M_S \cdot c_{bio} +v_{bio}'\cdot(c_{S,Feed}-K_S \frac{v_{bio}'}{v_{bio}^{max}-v_{bio}'})\\
		&&0&=v_{P}'\cdot M_{P}\cdot c_{bio} -c_{P}\cdot v_{bio}'\\
		&&\boldsymbol{v}' &\in \text{arg }\underset{\boldsymbol{v}\in \mathbb{R}^n}{\text{max}} 
		\quad v_{bio}   &&\\
		&&&  \begin{aligned} \text{s.t.} \quad \boldsymbol{S}\boldsymbol{v}&=\boldsymbol{0} \\
			v_{bio}&\geq f\cdot v_{bio,WT} \\
			\boldsymbol{v}&\geq \boldsymbol{v^{lower}}\circ (\boldsymbol{B}\boldsymbol{y})\\
			\boldsymbol{v}&\leq \boldsymbol{v^{upper}}\circ(\boldsymbol{B}\boldsymbol{y}),
		\end{aligned} 
	\end{aligned}
\end{equation} 
where $\circ$ denotes the element-wise product, $\bm{y}\in\{0,1\}^{r}$ denotes the binary knockout vector, and $r$ denotes the number of reactions in a metabolic network, which contains irreversible and reversible reactions, further referred to as a reversible network. 
The parameter $K$ is the number of maximal allowed knockouts, and the parameter $f$ denotes a value between zero and one. 
The matrix $\bm{B}$ maps the reactions in the irreversible metabolic network with those in the reversible network. 
The apostrophe $'$ denotes the cellular variables in the upper level. 
Note that, in the algorithm, they are fixed by the lower-level program.  

As described above, we chose to eliminate $D$ and $c_S$ from the formulation. 
In theory, the process variables $c_{bio}$ and $c_P$ could also be eliminated, as the equality constraints at the upper level are fixing them. 
Nevertheless, we decided on this formulation because preliminary tests and our intuition suggest that it results in simpler convergence. 
Note that our recent research indicates that an intermediate between full-space (keeping all variables) and reduced-space formulation (eliminating all variables except the degrees of freedom) is most beneficial in some cases \citep{Bongartz.2019, Najman.2021}. 
With $c_{bio}$ and $c_P$ being fixed by equality constraints, the substrate concentration in the feed $c_{S, Feed}$ is the only degree of freedom in the process conditions.
At the cellular level, $v_S$ will reach its upper bound based on an experimental value, as it is standard in other strain optimization examples (e.g., \citep{Burgard.2003}). 
Even if the Monod kinetics constrain $v_{bio}$ and $c_{S,Feed}$ is fixed during optimization; spare carbon atoms will go into product formation.  

In our study, the Monod constant was set to \SI{0.044}{\gram\per\liter}, and $v_{bio}^{max}$ was set to \SI{0.73}{\per\hour}, according to \citet{Wick.2001}.
The wild-type biomass flux $v_{bio,WT}$ was determined by a previously performed FBA, and the parameter $f$ was set to 0.1.
Furthermore, the ATP maintenance reaction flux threshold was set to \SI{6.86}{\milli\mol\per\gramdryweight\per\hour}, according to the default value in \textit{i}ML1515 \citep{Monk.2017}.
The upper bound of the glucose exchange reaction was set to \SI{10}{\milli\mol\per\gramdryweight\per\hour}.

Alternatively, we used Michaelis-Menten, the simplest form of enzyme kinetics.
Enzymes catalyze reactions within a cell, with metabolites as reactants. 
An enzyme-catalyzed reaction can be described using Michaelis-Menten kinetics. 
In alignment with \citet{Mahadevan.2002} and \citet{Ploch.2019}, we chose to apply Michaelis-Menten kinetics to the substrate uptake reaction.
Thus, the kinetics link the substrate concentration in the bioreactor with the substrate uptake rate.
More precisely, the substrate uptake rate is a function of the substrate concentration: $v_S = g(c_S)$.
The kinetics read
\begin{align*}
	v_{S} &= v_{S}^{max}\frac{c_S}{c_S+K_{S,MM}}. &(Michaelis-Menten)
\end{align*}
The maximal substrate uptake is denoted with $v_S^{max}$ in \si{\milli\mol\per\gramdryweight\per\hour}, and the Michaelis constant is $K_{S, MM}$ in \si{\gram\per\liter}. 
Again, we assumed that the kinetic parameters were constants and, thus, can be retrieved from the literature.
As $v_S$ is a cellular variable, the kinetics were implemented in the lower level. 
SimulKnock with embedded Michaelis-Menten kinetics reads 
\begin{equation}
	\label{eq:SimulKnock_MM}
	\begin{aligned}
		\underset{\boldsymbol{y},c_{S,Feed},c_S, c_{P}, c_{bio}, v_S', v_{bio}', v_P'}\max &&\quad c_{P}\cdot v_{bio}'  &&          & \\
		\text{s.t. }                              &&\sum_{i=1}^{r} (1-y_i)& \leq K&& \\
		&&0&=-v_{S}'\cdot M_S \cdot c_{bio} +v_{bio}'\cdot(c_{S,Feed}-c_S)\\
		&&0 &=v_{P}'\cdot M_{P}\cdot c_{bio} -c_{P}\cdot v_{bio}'\\
		&&\boldsymbol{v}' &\in \text{arg } \underset{\boldsymbol{v}\in \mathbb{R}^{n}}{\text{max}} 
		\quad v_{bio}   &&\\
		&&&\begin{aligned} \text{s.t.} \quad  \boldsymbol{S}\boldsymbol{v}&=\boldsymbol{0} \\
			v_{bio} &\geq f\cdot v_{bio,WT} \\
			v_{S} &= v_{S}^{max}\frac{c_S}{c_S+K_{S,MM}}\\
			\boldsymbol{v} &\geq \boldsymbol{v^{lower}}\circ(\boldsymbol{B}\boldsymbol{y})\\
			\boldsymbol{v}&\leq \boldsymbol{v^{upper}}\circ(\boldsymbol{B}\boldsymbol{y}).\\
		\end{aligned}
	\end{aligned}
\end{equation}
Again, $c_{bio}$ and $c_P$ are being fixed by the equality constraints in the upper level. 
However, unlike the formulation including Monod, $c_S$ is a degree of freedom at the process level now. It is included in the upper-level equality constraint in the summand $v_{bio}\cdot(c_{S,Feed}-c_S)$. With $v_{bio}$ being part of the upper-level objective, this term also has an optimum, and $c_S$ is determined during the optimization.
Another difference is that $v_S$ is now constrained by the kinetics instead of reaching the upper bound set by the user. 

In our study, the Michaelis constant was calculated with \SI{0.53}{\milli\mol\per\liter} multiplied by the molar weight of glucose, and $v_{S}^{max}$ was set to \SI{10}{\milli\mol\per\gramdryweight\per\hour}, both values according to \citet{Meadows.2010}.

To simplify the comparison of Monod and Michaelis-Menten, Table \ref{tab:comparison_Monod_MM} summarizes the different aspects of both kinetics applied in SimulKnock. 

\begin{table}[hbt]
	\centering
	\caption{Comparison of Monod and Michaelis-Menten in SimulKnock. Abbreviations: $c_S$: substrate concentration; $v_{bio}$: growth rate; $v_S$: substrate uptake rate; $f$, $g$: notation for functions.}
	\label{tab:comparison_Monod_MM}
	\resizebox{\textwidth}{!}{
		\begin{tabular}{lll}
			\hline
			&Monod                           &Michaelis-Menten\\
			\hline
			functionality    &$c_S = f(v_{bio})$              &$v_S = g(c_S)$  \\
			level    &whole-cell kinetics $\rightarrow$ process level & enzyme kinetics $\rightarrow$ cellular level \\
			$v_S$    &unrestricted, reaches its upper bound       &restricted by kinetics and $c_S$\\
			$c_S$    &calculated by kinetics    &determined by equality constraint and objective in process level \\
			assumption& gene deletion does not change growth behavior& gene deletion does not change substrate uptake behavior\\
			parameter availability &higher   &lower\\
			possible extensions& inclusion of inhibitions & consideration of multiple substrates\\
		\end{tabular}
	}
\end{table}

Apart from the aspects already discussed above, we found that, in general, the parameter availability for Monod is higher than for Michaelis-Menten. 
Furthermore, the inclusion of inhibitions is supposedly a straightforward extension with Monod, whereas including multiple substrates would potentially be easier with Michaelis-Menten (cf. \citep{Ploch.2019}).

\subsection{Reformulation of the bilevel program to a single-level quadratically-constrained quadratic program} 

We used the software packages libALE \citep{Djelassi.2019b} and libDIPS \citep{Jungen.2023} to implement SimulKnock in both versions---Monod \eqref{eq:SimulKnock_Monod} and Michaelis-Menten \eqref{eq:SimulKnock_MM}---and successfully solved the bilevel programs for the \textit{E. coli} core network \citep{Orth.2010} using the algorithm \citep{Mitsos.2008} without KKT-based tightening and with Gurobi v10.0 \citep{gurobi} as the subsolver.
However, this algorithm was too expensive for the genome-scale metabolic model \textit{i}ML1515. \citep{Monk.2017}. 
Consequently, we reformulated the bilevel programs to single-level mixed-integer quadratically-constrained quadratic programs (MIQCQP).
We conducted two reformulation steps for numerical reasons: first, we reformulated SimulKnock as a single-level program. 
Second, we eradicated nonlinear terms to achieve an MIQCQP formulation of SimulKnock. 

To the first reformulation: In both versions of SimulKnock (Monod or Michaelis-Menten embedded), the cellular level is linear (in the lower-level variables, as the process level fixes the process variables). 
Thus, we can apply strong duality and reformulate to a single-level program.  
The reformulations for SimulKnock with Monod kinetics embedded are identical to those presented in OptKnock \citep{Burgard.2003}, as the Monod kinetics are in the upper-level program. 
In the following, we will show the reformulation of SimulKnock with Michaelis-Menten kinetics embedded. 
For ease of notation, we introduced 
\begin{equation*}
	\boldsymbol{\tilde{C}}=\begin{bmatrix}
		\boldsymbol{I}^{|n|}\\
		\boldsymbol{I}^{|n|}\\
		\boldsymbol{c}^T_{kin}
	\end{bmatrix},\; \boldsymbol{\tilde{b}}=\begin{bmatrix}
		\boldsymbol{v^{lower}}\circ(\boldsymbol{B}\boldsymbol{y})\\-\boldsymbol{v^{upper}}\circ(\boldsymbol{B}\boldsymbol{y})\\d
	\end{bmatrix},\;c_i=\begin{cases} 0, \text{if}\; i\neq bio\\
		1, \text{if}\; i=bio,\end{cases} \forall i \in {1, ..., n},
\end{equation*}

where $\boldsymbol{c}$ describes the position of the biomass reaction and $\boldsymbol{I}$ is the unity matrix. The vector $\boldsymbol{c}_{kin}\in\{0,1\}^{|n|}$ has exactly one nonzero entry at the index of the glucose uptake reaction.
The scalar $d$ describes the right-hand side of the kinetics, i.e., $d = v_{S}^{max}\frac{c_S}{c_S+K_{S,MM}}$. 
By applying the definitions above, the lower-level program of \eqref{eq:SimulKnock_MM} can be rewritten as
\begin{equation}
	\label{eq:cellular_level_kinetics_canon}
	\begin{aligned}
		& \underset{\boldsymbol{v}}{\text{max}}
		& &   \boldsymbol{c}^T\boldsymbol{v}\\
		& \text{s.t.}
		& & \boldsymbol{S}\cdot \boldsymbol{v}=\boldsymbol{0}\\
		& & &\boldsymbol{\tilde{C}}\cdot \boldsymbol{v} \leq \boldsymbol{\tilde{b}}\\
		& & & \boldsymbol{v}\geq 0.
	\end{aligned}
\end{equation}

Converting program \ref{eq:cellular_level_kinetics_canon} to its dual and applying the strong duality theorem yields the system of equations
\begin{equation}
	\label{eq:cellular_level_dualized}
	\begin{aligned}
		\boldsymbol{c}^T\boldsymbol{v}&= \boldsymbol{\tilde{b}}^T\cdot \boldsymbol{\tilde{\mu}}\\
		\boldsymbol{S}\cdot\boldsymbol{v}&=\boldsymbol{0}\\
		\boldsymbol{\tilde{C}}\cdot \boldsymbol{v} &\leq \boldsymbol{\tilde{b}}\\
		\boldsymbol{c} &= \boldsymbol{S}^T\cdot \boldsymbol{\lambda} + \boldsymbol{\tilde{C}}^T\cdot\boldsymbol{\tilde{\mu}}\\
		&\boldsymbol{\tilde{\mu}}\geq \boldsymbol{0},
	\end{aligned}
\end{equation}
with $\boldsymbol{\lambda}\in\mathbb{R}^{|m|}$ being the dual variables corresponding to the mass balances for the metabolites and $\boldsymbol{\tilde{\mu}}\in\mathbb{R}_{+}^{|2n+1|}$ being the dual variables corresponding to the bounds set on the reaction fluxes through the inequality constraints and the kinetics.
The set of equations \ref{eq:cellular_level_dualized} can readily replace the cellular level in \eqref{eq:SimulKnock_MM}. 
The resulting program is a single-level mixed-integer nonlinear program. 
The nonlinearity is introduced by the fact that process and cellular variables are at the same level now.
The nonlinear term is the kinetic term in $d$.
These considerations about nonlinearity are also applicable to the reformulated SimulKnock with Monod embedded. 

In the second reformulation step, we reformulated the nonlinear kinetic term.
We introduced an additional variable for the fraction on the right-hand side of the kinetics, thereby introducing an additional equality constraint in the problem. 
For example, in the Michaelis-Menten kinetics, the fraction is replaced with the optimization variable $\sigma_{MM}$, and the constraints
\begin{equation}
	\begin{aligned}
		&v_{S}= v_S^{max} \sigma_{MM} \\
		&\sigma_{MM} (c_S+K_{S,MM})= c_S
	\end{aligned} 
	\label{eq:reform_MM}
\end{equation}
are introduced, which results in the single-level optimization problem becoming an MIQCQP. This reformulation is beneficial because MIQCQPs are solved directly by commercial solvers, such as Gurobi \citep{gurobi}.

Replacing the lower level in (\ref{eq:SimulKnock_MM}) with the set of constraints of (\ref{eq:cellular_level_dualized}) and applying the reformulation of the kinetic term suggested in (\ref{eq:reform_MM}) yields the single-level mixed-integer quadratic form of SimulKnock with Michaelis-Menten embedded:

\begin{equation}
	\label{eq:SimulKnock_MM_dualized}
	\begin{aligned}
		\underset{\boldsymbol{y},c_S,c_{S,Feed},\sigma}\max             &\quad c_{P}\cdot v_{bio}  &&          & \\
		\text{s.t.: } &\sum (1-y_i) \leq 1\\
		&0=-v_{S}\cdot M_S \cdot c_{bio} +v_{bio}\cdot(c_{S,Feed}-c_S)\\
		&0=v_{P}\cdot M_{P}\cdot c_{bio} -c_{P}\cdot v_{bio}\\
		&v_{S}= v_S^{max} \sigma_{MM}\\
		&c_S=\sigma_{MM}(c_S+K_{S,MM})\\
		&v_{bio} = \left[-\boldsymbol{v^{lower}}\circ\left(\boldsymbol{B}\boldsymbol{y}\right)\quad \boldsymbol{v^{upper}}\circ\left(\boldsymbol{B}\boldsymbol{y}\right)\quad v_S^{max} \sigma \right]\cdot \boldsymbol{\mu}\\
		&\boldsymbol{S}\cdot \boldsymbol{v}=\boldsymbol{0}\\
		&\boldsymbol{v}\geq \boldsymbol{v^{lower}}\circ\left(\boldsymbol{B}\boldsymbol{y}\right)\\
		&\boldsymbol{v}\leq \boldsymbol{v^{upper}}\circ\left(\boldsymbol{B}\boldsymbol{y}\right)\\
		&\boldsymbol{S}^T\cdot \boldsymbol{\lambda} +\left[\boldsymbol{I}^{|n|}\quad\boldsymbol{I}^{|n|}\quad \boldsymbol{c}_{S}\right] \cdot \boldsymbol{\mu} = \boldsymbol{c}\\
		&v_{bio}\geq f\cdot v_{bio}^{max}\\
		&\boldsymbol{v},\boldsymbol{\mu}\geq \boldsymbol{0}
	\end{aligned}
\end{equation}

Similar reformulations also reduce SimulKnock with Monod kinetics to an MIQCQP. 
The reformulated programs were implemented using Pyomo \citep{bynum2021pyomo, hart2011pyomo} and solved with Gurobi v10.0 \citep{gurobi}.

\section{Case studies}
We performed three case studies with SimulKnock. 
The first case study illustrates SimulKnock's mode of action compared to OptKnock and a sequential optimization approach. The second case study employs the genome-scale metabolic model \textit{i}ML1515 \citep{Monk.2017} of \textit{E. coli} for one to three knockout predictions and elaborates on the differences between Monod and Michaelis-Menten kinetics. 
The third case study compares the results of SimulKnock with the published results of an experimental study on \textit{E. coli} continuous fermentation \citep{vanHeerden.2013}.
The first case study was performed on up to 8 Intel Xeon E5-2640 CPU threads and was solved within less than a second. Case studies two and three were solved within five hours on the RWTH high-performance computing cluster using 48 threads and up to 4 GB of memory per thread.

\subsection{SimulKnock predicts different knockouts than OptKnock with embedded illustrative network}
\begin{figure}[htb]
	\centering
	\includegraphics[trim = 6cm 4cm 6.0cm
	1.5cm, clip, width=\textwidth]{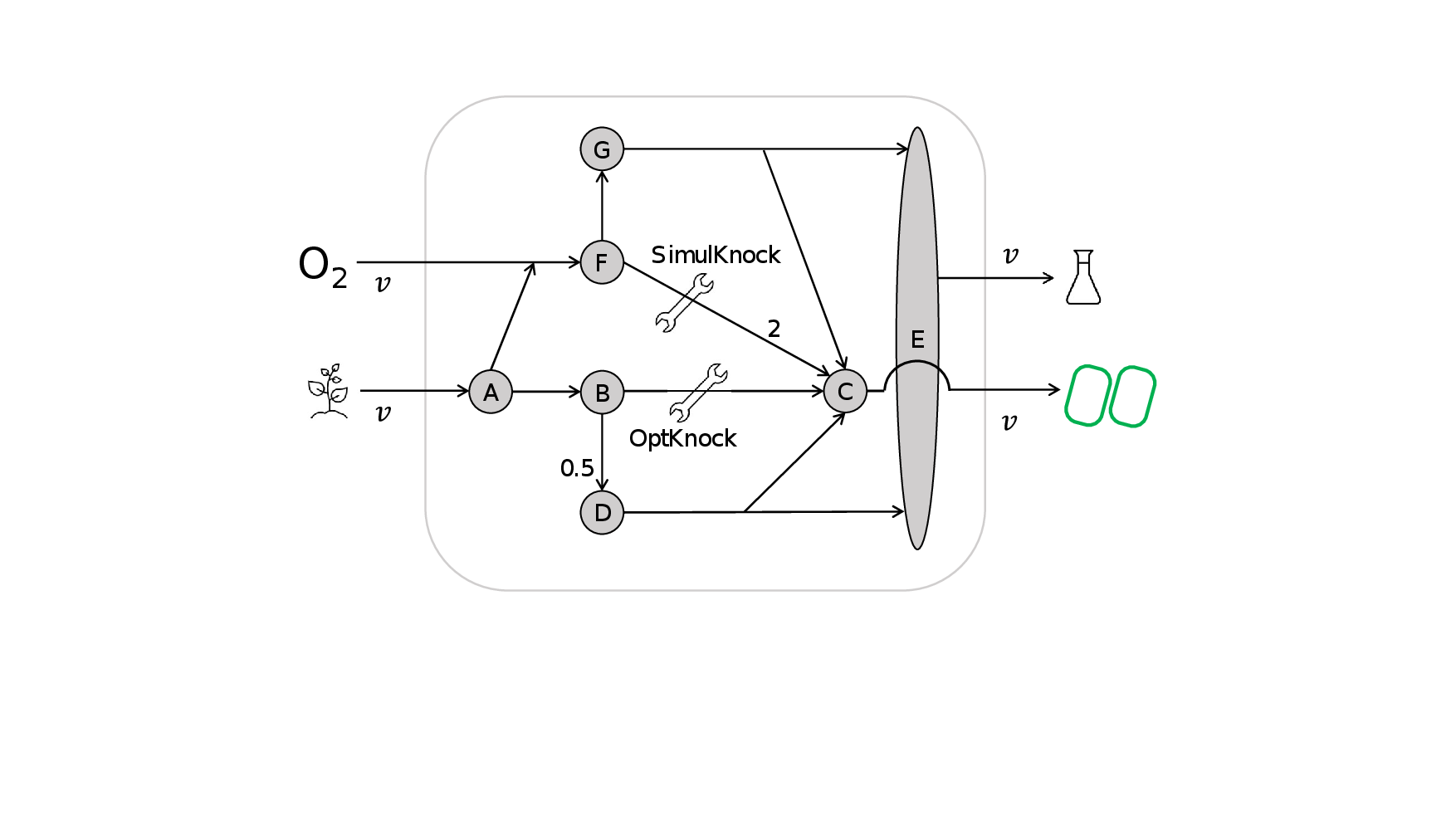}
	\caption{Illustrative network. The grey circles A to G denote metabolites, $v_O$ and $v_S$ denote the flux of oxygen and the substrate uptake reaction, respectively; $v_{P}$ and $v_{bio}$ are the target chemical and biomass flux, respectively. All stoichiometric coefficients are one, except where indicated with numbers. For one maximum allowable knockout, the wrenches indicate the knockout prediction of SimulKnock and OptKnock.}
	\label{fig:example_network}
\end{figure}

To show SimulKnock's mode of action compared to sequential optimization using OptKnock, we constructed an illustrative network (Figure~\ref{fig:example_network}). 

We performed four optimizations with the illustrative network embedded: $(i)$ an FBA with maximization of biomass as the objective (referred to as wild-type), $(ii)$ an OptKnock prediction, $(iii)$ a sequential approach using OptKnock (see below), and $(iv)$ SimulKnock.
Table~\ref{tab:results_toy_network} shows the results of the case study. 

\begin{table}[ht]
	\centering
	\caption{Results of the case study with the illustrative network (Figure~\ref{fig:example_network}) embedded. The substrate is glucose, and the maximum substrate flux allowed is set to $v_S^{upper} =\SI{10}{\milli\mol\per\gramdryweight\per\hour}$; the upper bounds on A-B and A-F are set to \SI{7}{\milli\mol\per\gramdryweight\per\hour} and \SI{3}{\milli\mol\per\gramdryweight\per\hour}, respectively; the maximum oxygen flux allowed is $v_O^{upper} =\SI{3}{\milli\mol\per\gramdryweight\per\hour}$. Michaelis-Menten kinetics were applied in SimulKnock and the sequential approach, with the Michaelis constant $K_{S,MM}=\SI{0.53}{\milli\mol\per\liter}$ \citep{Meadows.2010}. In the sequential approach, 1. OptKnock, and 2. A process optimization is performed with embedded knockouts from OptKnock. Abbreviations: $v_S$: substrate uptake flux, $v_O$: oxygen flux, $v_{bio}$: growth rate, $c_{bio}$: biomass concentration, $c_P$: product concentration, STY: space-time yield}
	\label{tab:results_toy_network}
	\begin{tabular}{lllllllll}
		\toprule
		&Knockout & $v_S$ & $v_O$ & $v_{bio}$ & $v_{P}$ &$c_{bio}$ & $c_{P}$ &STY\\
		& [-] & \multicolumn{4}{c}{[\si{\milli\mol\per\gramdryweight\per\hour}]} & \multicolumn{2}{c}{[\si{\gram\per\liter}]} & [\si{\gram\per\liter\per\hour}]\\
		\midrule
		Wild-type & -& 10&3&13&0&-&-&-\\
		OptKnock &B-C&10&3&9.5&3.5&-&-&-\\
		Sequential &B-C&7.8&3&8.4&2.4&8.7&2.5&21.0 \\
		SimulKnock &F-C&3&3&3&3&9.8&9.8&29.3\\
		\bottomrule	
	\end{tabular}
\end{table}

The sequential approach denotes a two-step procedure. 
First, OptKnock is performed. 
Second, the predicted knockouts are applied to the network, and a continuous process optimization is performed. 
The continuous process optimization looks like the formulation of SimulKnock but with fixed knockouts, which were determined from OptKnock. 
Thus, OptKnock is included in the sequential approach. 
Still, we added the results of OptKnock for completeness to allow for a direct comparison of functionalities with SimulKnock. 
Michaelis-Menten kinetics were applied in the sequential approach and SimulKnock. 
The FBA and OptKnock predictions do not imply process variables, e.g., the substrate or product concentration, which is also why the kinetics were not applied in OptKnock or FBA.

The results indicate that SimulKnock and OptKnock predict different optimal knockouts for one maximum allowable knockout. 
The reason for this difference is that OptKnock optimizes for the target chemical production, whereas SimulKnock optimizes for space-time yield. 
Indeed, the product flux $v_{P}$ is higher for OptKnock than for SimulKnock. 
The applied kinetics in SimulKnock and the sequential approach force the substrate uptake lower than the maximum allowed value. 
Namely, in the sequential approach, the substrate uptake is set to \SI{8}{\milli\mol\per\gramdryweight\per\hour} by the process optimization, whereas OptKnock predicted it to be \SI{10}{\milli\mol\per\gramdryweight\per\hour}.
In SimulKnock, only the reactions in the upper part of the network, i.e., the reactions going via the metabolites A, F, G, E, and C, are active and used for biomass and target chemical production. 
In contrast, in OptKnock and the sequential approach, the reactions going via B, D, and E are used for target chemical production, while the route via F and C produces additional biomass. 
The space-time yield, which is the product of the growth rate $v_{bio}$ and the product concentration $c_{P}$, is higher for SimulKnock than for the sequential approach.
Note that both formulations predict the same two knockouts, B-C and F-C, and the same space-time yield for two maximum allowable knockouts. 
However, this does not generally indicate that with more allowed knockouts, the results of OptKnock and SimulKnock would be similar. 
Instead, the illustrative network is constructed so small that no other options exist. 

\subsection{SimulKnock can achieve higher space-time yields than sequential optimization with embedded \textit{i}ML1515}

To further investigate the computational tractability of the SimulKnock approach, we chose a genome-scale metabolic model, \textit{i}ML1515 \citep{Monk.2017}, which describes the metabolism of \textit{E. coli}. 
The metabolic network includes 1516 genes, resulting in 1877 metabolites and 2712 metabolic reactions. 
The objective was set to maximize the space-time yield of six target chemicals for a glucose-limited medium.\\
Furthermore, to highlight the superiority of the simultaneous strain and process optimization, SimulKnock was compared against the sequential approach on the space-time yield of the target chemical.\\

\subsubsection{One knockout}
For one knockout prediction, SimulKnock and the sequential approach furnished identical results. 
The knockout predictions and space-time yields are in the Supporting Information (Figure~\ref{fig:comparison_1KO} and Table~\ref{tab:data_comparison_1KO}).

\subsubsection{Two knockouts: Monod vs. Michaelis Menten}
\begin{figure}[htb]
	\centering
	\begin{subfigure}[b]{0.49\textwidth}
		\centering
		\includegraphics[width=\linewidth]{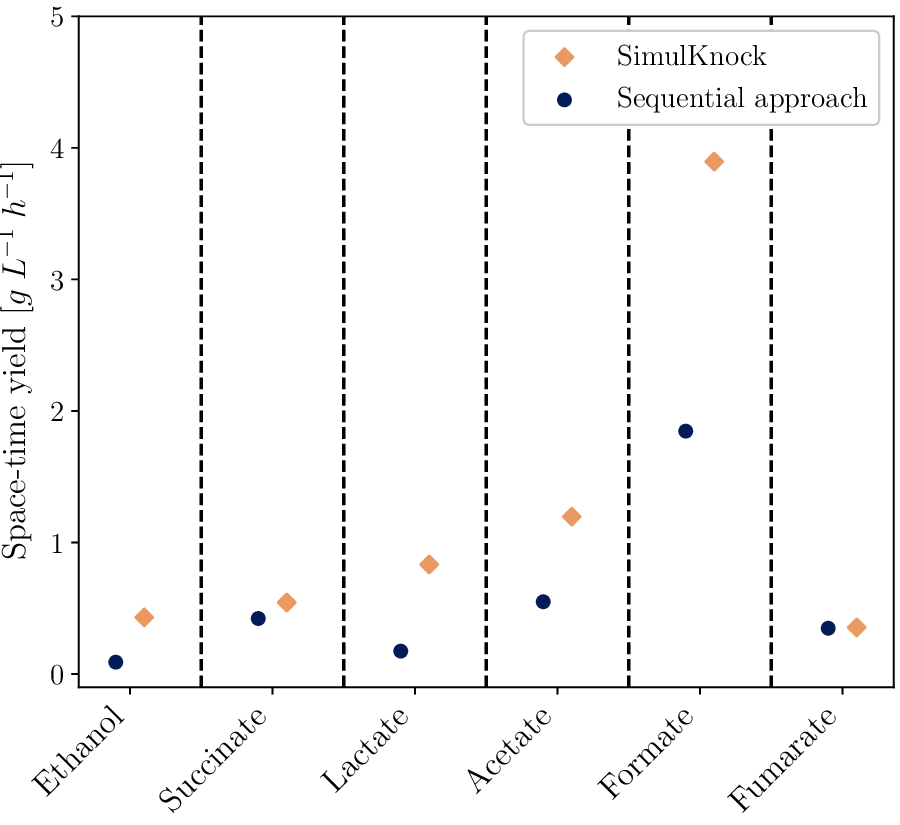}
		\caption{$Michaelis-Menten$}
		\label{subfig:comparison_2_MM}
	\end{subfigure}
	\hfill
	\begin{subfigure}[b]{0.49\textwidth}
		\centering
		\includegraphics[width=\linewidth]{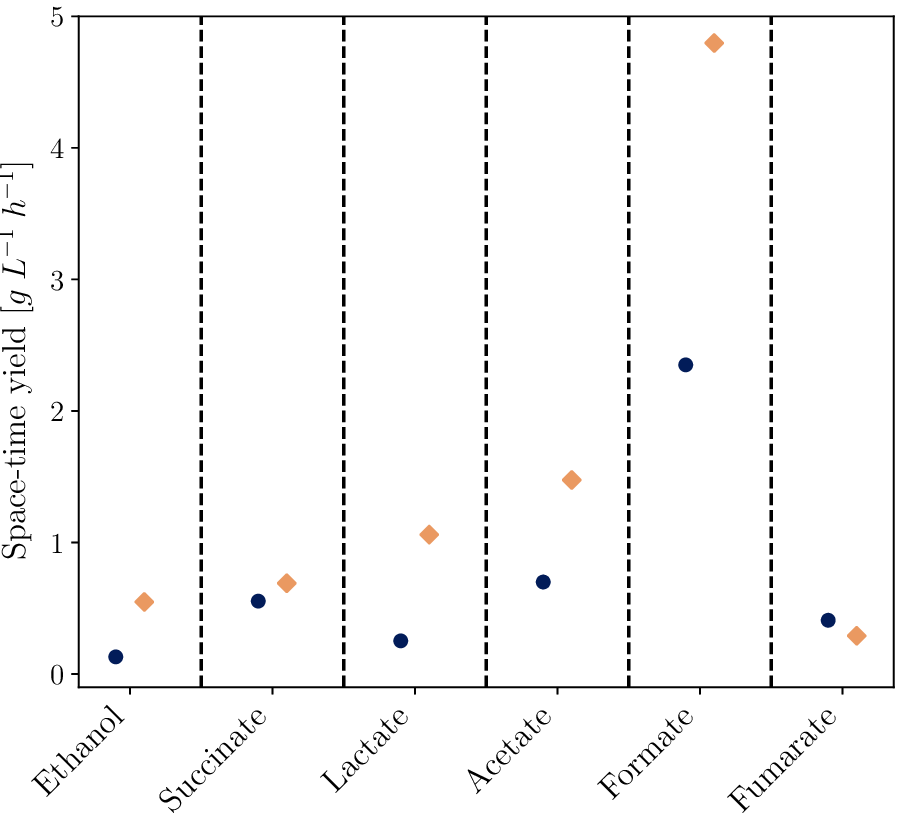}
		\caption{$Monod$}
		\label{subfig:comparison_2_Monod}
	\end{subfigure}
	\hfill
	\caption{Comparison of space-time yield using Michaelis-Menten and Monod kinetics for two reaction eliminations on the genome-scale metabolic model \textit{i}ML1515 \citep{Monk.2017}. Ethanol, succinate, and lactate were produced via anaerobic pathways; acetate, formate, and fumarate were produced via aerobic pathways.}
	\label{fig:comparison_2}
\end{figure}

In this case study, we allowed for two gene knockouts.
The comparison results of Michaelis-Menten and Monod kinetics are plotted in Figure~\ref{fig:comparison_2}.

Note that for \textit{E. coli}, the growth rate is higher with oxygen than without. The lower-level optimization problem always activates the oxygen uptake to maximize growth. 
Hence, in SimulKnock, the oxygen supply must be by the modeler. 
We performed runs with and without an oxygen supply. 
Figure~\ref{fig:comparison_2} shows the results where SimulKnock achieved the higher absolute value.  

SimulKnock predicts significantly higher space-time yields for four of the six considered target chemicals with embedded Michaelis-Menten kinetics.
The difference in the space-time yield from the two approaches is due to different reaction knockouts predicted (see the knockouts table \ref{tab:data_comparison_2KO} in the Supporting Information for reference). 
For example, SimulKnock targets the 6-phosphogluconolactonase and the acetate reversible transport for higher ethanol space-time yield instead of the ATP synthase and the triose-phosphate isomerase.
The most significant difference becomes visible for formate. 
Instead of the phosphoglycerate mutase, SimulKnock predicts the phosphofructokinase will be knocked out. 
Even if SimulKnock and the sequential approach do not predict the same knockouts for fumarate, they reach the same result with Michaelis-Menten. 

The kinetics also influence the knockout predictions. 
SimulKnock with embedded Michaelis-Menten kinetics suggests blocking Acetaldehyde dehydrogenase and D-alanine-D-alanine dipeptidase for succinate production. 
With Monod embedded dytosine deaminase is knocked out instead of the dipeptidase,. 
Interestingly, SimulKnock achieves a lower space-time yield for fumarate than the sequential approach with embedded Monod kinetics. 

The sequential formulation maximizes the target chemical concentration at the process level, with the biomass growth rate being used only in the metabolism level problem. Therefore, the targeted reaction knockouts decrease the growth rate to achieve the maximum target chemical concentration. 
In contrast, the SimulKnock formulation recognizes the trade-off between a higher growth rate and a higher target chemical flux. 
It is, therefore, not surprising that the sequential approach predicts a higher target chemical concentration and lower biomass concentrations than the SimulKnock approach.
The different modes of action become even more visible when comparing the molar product yield on the substrate. 
We computed the yield by dividing the target chemical flux $v_P$ by the substrate flux $v_S$ (cf. Table \ref{tab:data_comparison_2KO} in the Supporting Information for the data).
In all cases, the molar product yield was lower with SimulKnock than with the sequential approach. 
At maximum, in the case of acetate production, the yield was \SI{45}{\%} lower, whereas the space-time yield was 2.6 times higher.

\subsubsection{Three knockouts}
\begin{figure}
	\centering
	\includegraphics[width=0.49\linewidth]{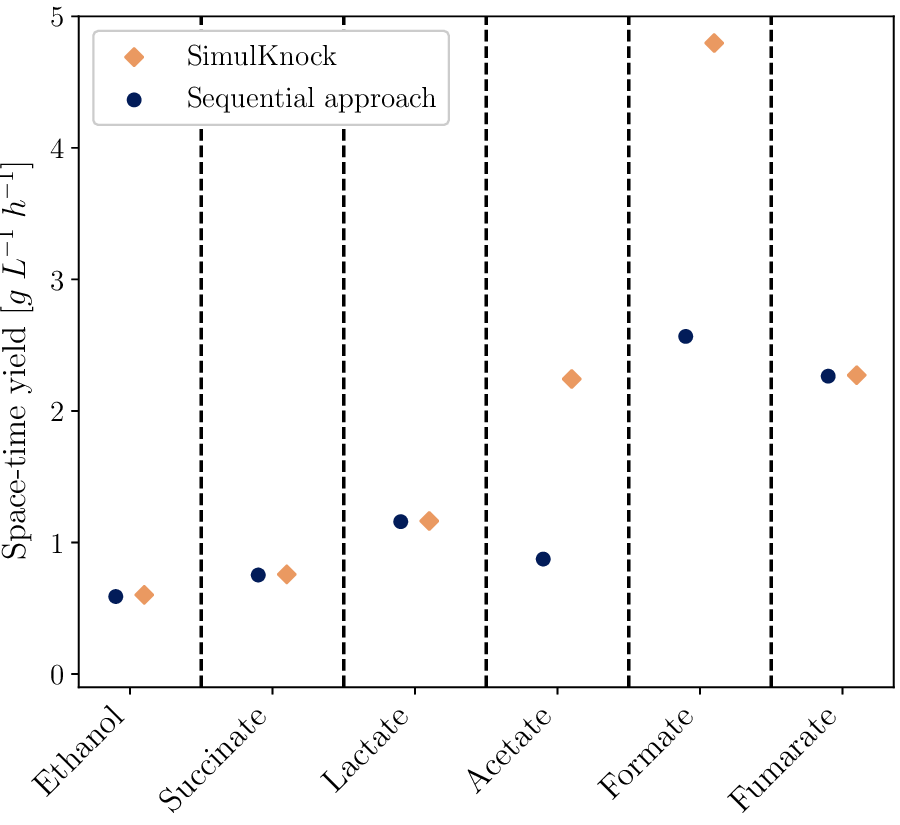}
	\caption{Space-time yield using Monod kinetics for three reaction eliminations on the genome-scale metabolic model \textit{i}ML1515 \citet{Monk.2017}. All six target chemicals were produced via aerobic pathways.}
	\label{fig:comparison_3}
\end{figure}

Upon increasing the possible reaction eliminations to three, the SimulKnock formulation with Michaelis-Menten kinetics could not be solved within a feasible time limit. 
Thus, only the Monod kinetics results are presented in Figure~\ref{fig:comparison_3}.
Note that, again, we performed runs with and without an oxygen supply. 
Figure~\ref{fig:comparison_3} shows the results where SimulKnock achieved the higher absolute value, all achieved using aerobic pathways. 

Acetate and formate space-time yield is increased by about 100\% when using SimulKnock, whereas ethanol, succinate, and fumarate increase marginally, and lactate remains the same. 
Once again, the reason for this difference is the different reaction knockouts predicted by the two approaches (cf. Table \ref{tab:data_comparison_3KO} in the Supporting Information). 
The knockouts for lactate and fumarate are similar, often targeting neighboring reactions, thereby predicting similar yields. 
When ethanol is considered, two of the three reaction eliminations target the same reactions with SimulKnock when compared to the sequential approach. 
When the target chemical is acetate, all three knockouts are different; with formate, two of the three knockouts are different.  
This observation hints towards the possibility of missing possible knockout strategies with the sequential approach.

It is worth noting that the assumption made in this study about the biomass reaction being the rate-limiting step has yet to be experimentally validated. 
Therefore, reaction eliminations could, theoretically, target the rate-limiting step on which the Monod kinetic parameters were fitted. 
In that case, a new parameter fitting using the modified organism would be required. In contrast, the Michaelis-Menten kinetics do not require parameter refitting due to reaction eliminations, as the substrate uptake reaction is always active.

\begin{figure}[htb]
	\centering
	\begin{subfigure}[b]{0.49\textwidth}
		\centering
		\includegraphics[width=\linewidth]{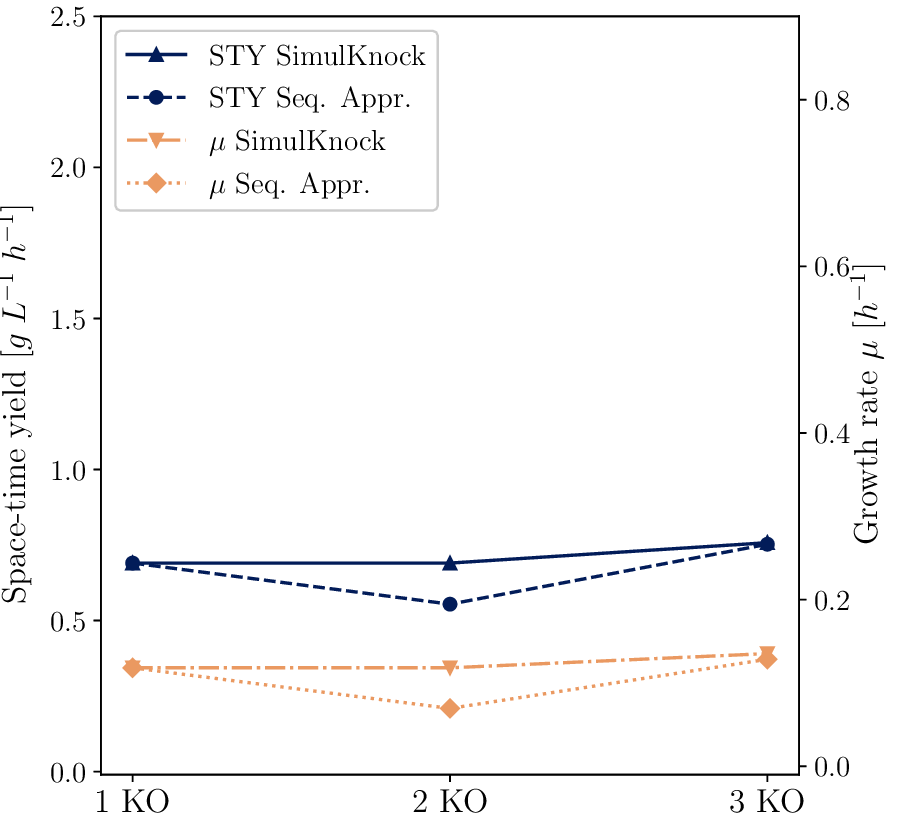}
		\caption{$Succinate$}
		\label{subfig:seq_monod_succ}
	\end{subfigure}
	\hfill
	\begin{subfigure}[b]{0.49\textwidth}
		\centering
		\includegraphics[width=\linewidth]{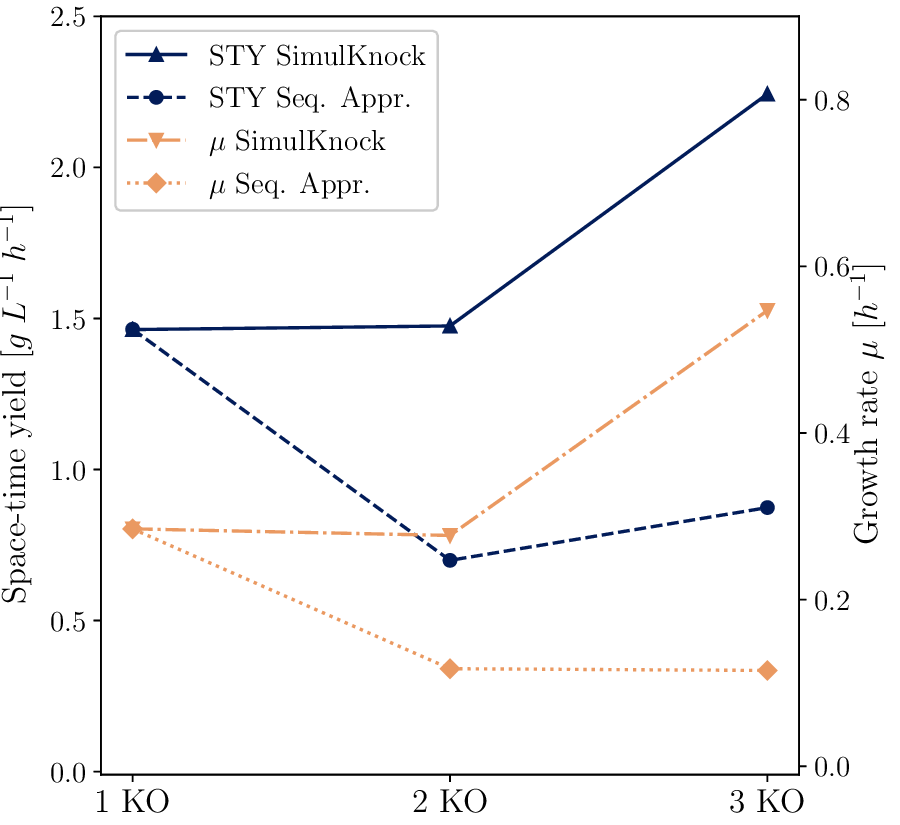}
		\caption{$Acetate$}
		\label{subfig:seq_monod_ac}
	\end{subfigure}
	\hfill
	\caption{Change of space-time yield and growth rate over the number of knockouts for \textit{i}ML1515 \citep{Monk.2017} with embedded Monod kinetics. Succinate was produced anaerobically for 1 and 2 knockouts and aerobically for 3 knockouts. Acetate was produced via an aerobic pathway. STY: space-time yield, Seq. Appr.: Sequential approach.}
	\label{fig:sty_growth_plot}
\end{figure}

When we line up space-time yield and growth rate with an increasing number of knockouts, we observe different behaviors: 
For ethanol, succinate, lactate, and formate, space-time yield and growth rate with SimulKnock stay more or less constant over one to three knockouts. 
Exemplary, Figure \ref{subfig:seq_monod_succ} shows the course for succinate. 

In the sequential approach, space-time yield decreases with an increasing number of knockouts, as depicted at two gene knockouts in Figure \ref{subfig:seq_monod_succ}. 
This decrease in space-time yield comes due to the mathematical structure of OptKnock. 
In OptKnock, high product yields are achieved at the expense of reduced growth (cf. Table \ref{tab:data_comparison_3KO} in the Supporting Information for data). 
In turn, space-time yield (equals product concentration times growth rate) is linearly dependent on the growth rate. 
The figure depicts the aerobic pathway for three gene knockouts, following our decision to show the results for the highest SimulKnock space-time yield. 
Hence, the fermentation conditions change compared to the anaerobic conditions that apply for one and two gene knockouts. 
While the space-time yield of the sequential approach would decline for anaerobic conditions (not shown in the figure), it reaches similar results as SimulKnock for aerobic conditions. 

Acetate and fumarate exhibit a different behavior, as depicted in an example in Figure \ref{subfig:seq_monod_ac}.
With SimulKnock, the space-time yield rises significantly at three knockouts, in alignment with the increase in growth. 
With the sequential approach, the space-time yield is lowest for two knockouts and then increases again, decreasing the growth rate. 
This interesting behavior stems from the two-step procedure. 
Only in the second step are the process conditions optimized, but they can only influence one of the two factors, namely the product concentration $c_P$. 
Thus, there is no monotonous behavior with the sequential approach. 
This finding displays the problems often occurring during scale-up when sequentially optimizing the microbial strain and the process conditions. 
On the other hand, acetate case underpins the potential of SimulKnock.

\subsection{SimulKnock can furnish meaningful knockout predictions but also exhibits model-experiment mismatch}

\begin{figure}[htb]
	\centering
	\includegraphics[trim = 0cm 0cm 0cmm 1.1cm, clip, width=\textwidth]{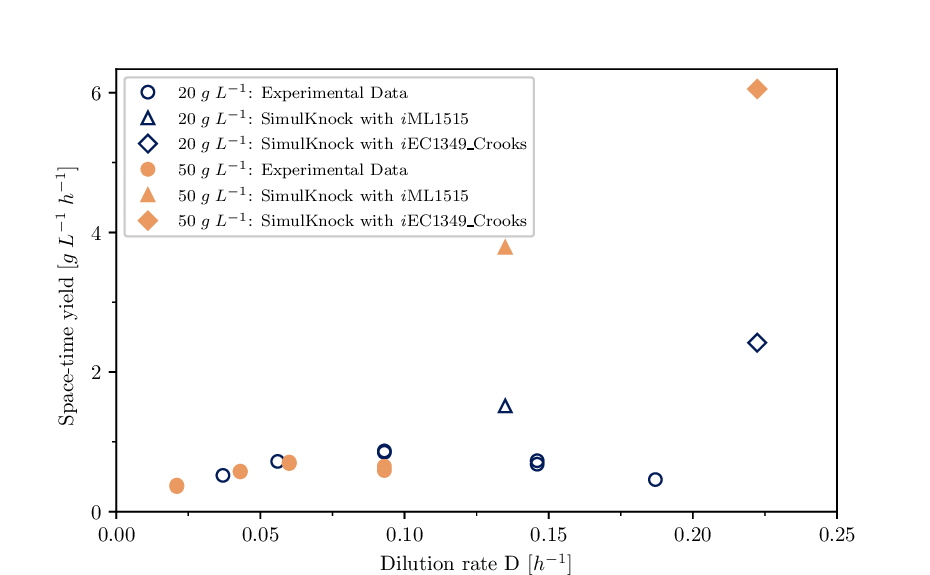}
	\caption{Comparison of experimental laboratory data for succinic acid production using \textit{E. coli} KJ134 \citep{vanHeerden.2013} with results from SimulKnock, based on the metabolic networks \textit{i}ML1515 \citep{Monk.2017} and \textit{i}EC1349\_Crooks \citep{Monk.2016}. Two glucose feed concentrations were studied: \SI{20}{\gram\per\liter} and \SI{50}{\gram\per\liter}.}
	\label{fig:Lab_Comparison}
\end{figure}

In the last case study, we took the experimental study of \citet{vanHeerden.2013} as a test case to elaborate on how far SimulKnock reproduces their results regarding space-time yield, dilution rate, and knockouts. 
We aimed to see whether SimulKnock furnishes experimentally meaningful results. 
In their experimental study, \citet{vanHeerden.2013} produced succinic acid from \textit{E. coli} KJ134 in a continuous fermentation with two different glucose feed concentrations.
We applied SimulKnock with embedded Monod kinetics to the \textit{E. coli} GEMs \textit{i}ML1515 \citep{Monk.2017} as well as \textit{i}EC1349\_Crooks \citep{Monk.2016} with three maximum allowable knockouts.
While \citet{vanHeerden.2013} tested different dilution rates and measured the space-time yield for each, SimulKnock predicted the optimal space-time yield and the corresponding dilution rate. 
Figure~\ref{fig:Lab_Comparison} shows the case study results, depicted as the space-time yield over the dilution rate.
Note that, in continuous fermentation, the dilution rate equals the growth rate $\mu$ and the biomass flux $v_{bio}$.
In the experimental study, 13 genetic modifications were applied to the organism to reach a maximum space-time yield of around \SI{1}{\gram\per\liter\per\hour} on \SI{20}{\gram\per\liter} glucose feed.
SimulKnock predicted a maximum space-time yield of about 4 and \SI{6}{\gram\per\liter\per\hour} with  \textit{i}ML1515 and \textit{i}EC1349\_Crooks, respectively, on \SI{50}{\gram\per\liter}. 
The space-time yields were reached with three knockouts.
The reactions to be eliminated, predicted with \textit{i}ML1515 embedded, were acetate kinase, ATP synthase, and fumarase.
With \textit{i}EC1349\_Crooks embedded, the predicted reaction eliminations were phosphotransacetylase, ATP synthase, and succinyl-CoA synthetase.
Amongst them, acetate kinase and phosphotransacetylase were also eliminated in the experimental strain. 
Hence, SimulKnock proves to predict experimentally meaningful results.
However, the figure reveals a substantial model-experiment mismatch. 
Especially the results of \textit{i}EC1349\_Crooks are higher than any experimental value. 
In \textit{i}EC1349\_Crooks, the high growth rate might result from the biomass function difference.
For illustration, for an FBA with a glucose limitation of \SI{10}{\milli\mol\per\gramdryweight}, \textit{i}ML1515 has a growth rate of \SI{0.87}{\per\hour}, whereas \textit{i}EC1349\_Crooks has a growth rate of \SI{1.02}{\per\hour}.
The mismatch in space-time yield that comes visible with both metabolic networks might result from SimulKnock not considering inhibition and repression effects and being overly optimistic by not considering byproduct formation unless the byproduct formation hinders biomass formation. 
The mismatch could be reduced by including experimental flux data of the organism in question, for example, from \citet{Fischer.2004} for \textit{E. coli}, or by using additional information from a protein allocation model \citep{Alter.2021}.
In both cases, the bounds of specific fluxes would be fixed in a preprocessing step, which can readily be done in the SimulKnock code. 
These refinements of the metabolic network and adaptations at the cellular level are interesting in case of application but were not in the scope of this work.

\section{Conclusion}
We presented SimulKnock, an optimization formulation combining continuous fermentation optimization with strain optimization via gene knockouts. 
Such an optimization formulation allows us to consider industrial fermentation conditions already in the strain design, thereby overcoming current scale-up difficulties. 
SimulKnock is a bilevel optimization formulation, which we transformed into a single-level mixed-integer nonlinear optimization program and solved globally for an illustrative example network as well as for the \textit{E. coli}  GEM \textit{i}ML1515 \citep{Monk.2017}.
We applied Monod or Michaelis-Menten kinetics to connect the process level to the cellular level in SimulKnock. 
SimulKnock with applied Monod kinetics showed similar results as Michaelis-Menten kinetics, with Monod being less computationally intensive.
SimulKnock predicted different knockouts than OptKnock \citep{Burgard.2003}. 
Also, the space-time yield was significantly higher with SimulKnock compared to a sequential approach of OptKnock plus process optimization for aerobic and anaerobic conditions. 
Compared to experimental data \citep{vanHeerden.2013}, SimulKnock indicated that higher space-time yields could be achieved with fewer knockouts.

The computations for SimulKnock with embedded GEM were performed on a high-performance computing cluster. 
A decrease in the computation demand and reduction of computation time could be achieved by applying further network reduction methods \citep{Erdrich.2015} or by considering essential genes \citep{Goodall.2018}.
Future work should consider that kinetics would have to be adjusted when knockouts are applied to the organism because Michaelis-Menten, describing the substrate uptake, and Monod, describing the growth behavior, are affected by genetic modifications. 
Especially for altered growth behavior, the update of the kinetics should be considered.
Our case study with experimental data indicated that SimulKnock furnishes overly optimistic results. 
The conversion to a robust optimization formulation, as suggested in RobustKnock \citep{Tepper.2010}, could tackle this observation. 
Moreover, experimental flux data (for \textit{E. coli}, see \citet{Pharkya.2006} and \citet{Fischer.2004}) could reduce the model-experiment mismatch by fixing the bounds of specific fluxes in a preprocessing step. 
To make a more accurate prediction, SimulKnock could be extended both at the process level and at the cellular level, e.g., by more elaborate microbial formulations \citep{Brochado.2012, Segre.2002, Landon.2019}, as depicted by OptKnock \citep{Burgard.2003} and its extensions \citep{Apaydin.2017, Kim.2010}.
These formulations include a reference flux distribution, which could be furnished using a protein allocation model \citep{Alter.2021} or experimental \textit{in vivo} flux data \citep{Kuepfer.2005}.
The extension could also comprise other genetic modifications, i.e., gene insertion \citep{Pharkya.2004} and regulation \citep{Pharkya.2006}, as well as other process operation modes, e.g., continuous mode with cell retention, batch, and fed-batch.

\section*{Author contributions}
ALZ designed the SimulKnock formulation, the reformulation, the implementation, and the case studies.
MS performed the formulation setup, the reformulation, the implementation, and the case studies under the supervision of ALZ and AsMa. 
ALZ and MS visualized the data. 
ALZ, MS, and AsMa analyzed and discussed the data. 
ALZ and AsMa wrote the manuscript draft. 
AlMi and LMB had the initial idea of SimulKnock, discussed the data, and reviewed the draft.
AlMi secured funding.
All authors read and approved the final manuscript.

CRediT: \\
ALZ: Conceptualization, Methodology, Software, Validation, 
Investigation, Writing - Original Draft (Lead), Writing - Review and 
Editing,  Visualization\\
AsMa: Methodology, Software, Investigation, Writing - Review and Editing\\
MS: Methodology, Software, Validation, Investigation, Visualization, 
Writing - Review and Editing\\
LMB: Conceptualization, Writing - Review and Editing\\
AlMi: Conceptualization, Supervision, Methodology, Writing- Review and 
Editing, Funding acquisition

\section*{Acknowledgments}
This project was funded by the Deutsche Forschungsgemeinschaft (DFG, German Research Foundation) under Germany´s Excellence Strategy – Cluster of Excellence 2186 ``The Fuel Science Center'' – ID: 390919832.
This project benefited from the work of Clemens Kortmann, who improved the modularity and the user-friendliness of our git repository. 
Computations were performed with computing resources granted by RWTH Aachen University under project thes1376.

\section*{Data Availability and Reproducibility Statement}
The implementation of the SimulKnock approach, the data pre-processing, and the interface with the solver Gurobi are openly available in our GitLab repository ``SimulKnock: Simultaneous design of fermentation and microbe'' at \url{https://git.rwth-aachen.de/avt-svt/public/simulknock}. 
Data from the manuscript's figures are tabulated in the Supporting Information.

\section*{Declaration of conflicts}
No conflicts to declare.

\section*{Graphical Abstract}

	\begin{figure}[htb]
		\centering
		\includegraphics[width=\textwidth]{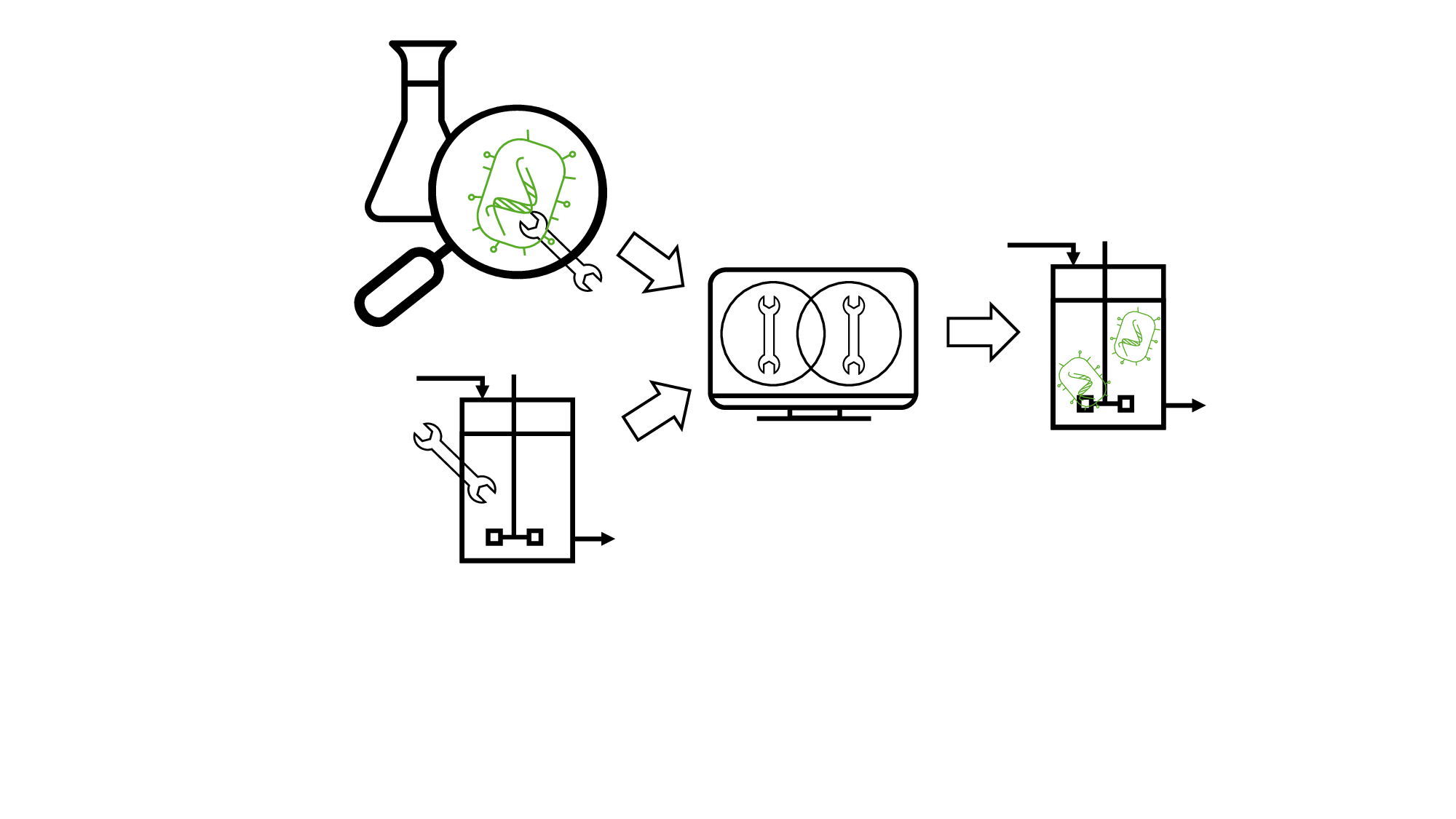}
		\caption{Microbial network and process optimizations are combined in our optimization formulation SimulKnock.}
	\end{figure}

\bibliographystyle{vancouver-authoryear.bst}
\bibliography{bibliography.bib}

\newpage
\setcounter{page}{1}
\appendix
\setcounter{table}{0}
\renewcommand{\thetable}{A.\arabic{table}}
\setcounter{figure}{0}
\renewcommand{\thefigure}{A.\arabic{figure}}
\setcounter{equation}{0}
\renewcommand{\theequation}{A.\arabic{equation}}

\section{Supporting Information}
\subsection{Analysis of the formulation with both kinetics embedded}
\label{seq:inclusion_both_kinetics}
Including both kinetics, Michaelis-Menten and Monod, may be tempting to tackle the drawback of Monod, which is that the substrate uptake flux $v_S$ is unrestricted. 
Regarding the parameter availability in the literature, the inclusion of two kinetics is impractical, as it typically is difficult to attain the parameters for both kinetics. 
Moreover, including both kinetics may result in an inconsistency, which we analyze in the following.

Looking at the mathematical functionalities, the inclusion of both kinetics leads to
\begin{equation*}
	\begin{aligned}
		&\text{Monod:\,} c_S = f(v_{bio}) \text{\,and Michaelis-Menten:\,} v_S = g(c_S) \\
		&\xRightarrow[]{} v_S = g(f(v_{bio}))\\
		&\text{metabolic\,network:\,} v_{bio} \in \mathcal{H}(v_S),
	\end{aligned}
\end{equation*}
where the set $\mathcal{H}$ indicates that $v_{bio}$ depends on $v_S$ through their connections in the metabolic network and that there is no unique solution. 
These relations could lead to inconsistency.
Alternatively, the kinetics could constrain $v_S$ and $v_{bio}$ to such an extent that the metabolic network would only simulate the product and side-product fluxes instead of optimizing $v_{bio}$. 
This case is not desirable either because we assumed that the metabolic network is more accurate than the kinetics, which is why we developed SimulKnock in the first place. 

Based on these considerations, we tested a formulation where we inserted Monod into Michaelis-Menten kinetics such that, as described above, $v_S = g(f(v_{bio}))$.
We conducted exemplary runs using \textit{i}ML1515.
The formulation was infeasible for anaerobic production, exemplarily shown for succinate (1 and 2 knockouts) and acetate (1 knockout).
The formulation is feasible only when ignoring the maintenance threshold, with the resulting space-time yield being zero. 
For aerobic production, the formulation did not converge within the time limit of \SI{20.000}{\second} for formate production. 
For succinate (2 knockouts) and acetate production (1 knockout), the achieved space-time yields were far lower (\SI{0.16}{\gram\per\liter\per\hour} and \SI{0.013}{\gram\per\liter\per\hour}, respectively) than the results with only one kinetics embedded.
At the same time, we observed elevated growth rates. 
For succinate production (2 knockouts), for example, the growth rate was \SI{0.68}{\per\hour}, which, set in relation to the experimental results of \citet{vanHeerden.2013}, is not realistic. 

From our mathematical analysis and the exemplary study, which aligned with the mathematical analysis, we conclude that it is not beneficial to include both kinetics.

\subsection{Data for one reaction elimination}
\begin{figure}[htb]
	\centering
	\begin{subfigure}[b]{0.49\textwidth}
		\centering
		\includegraphics[width=\linewidth]{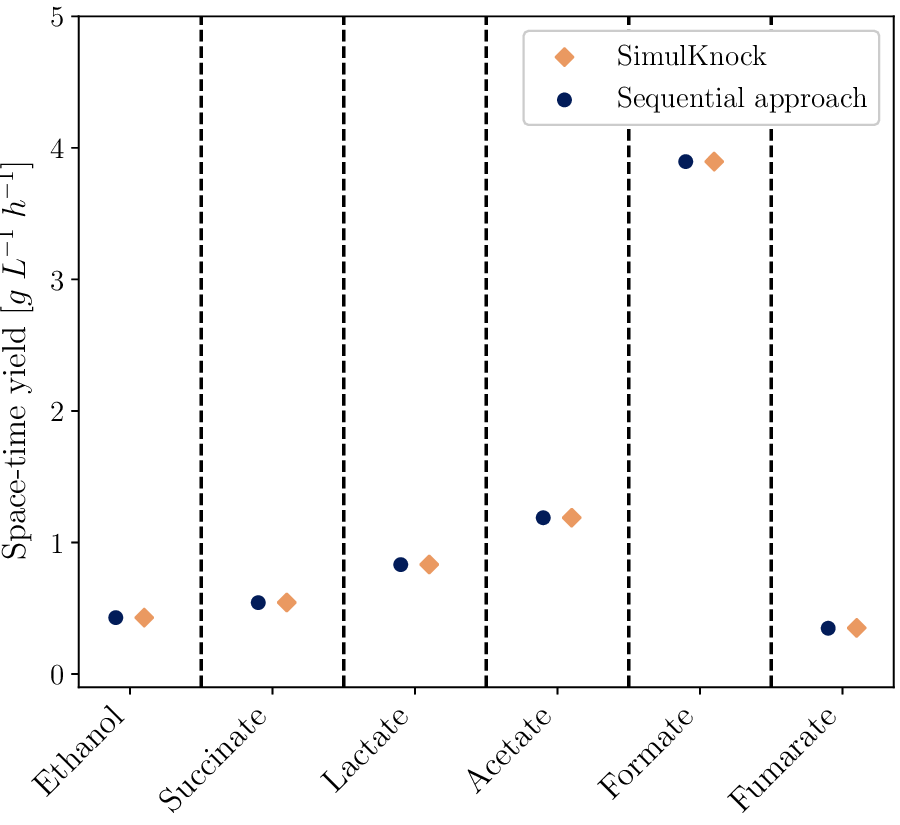}
		\caption{$Michaelis-Menten$}
		\label{subfig:comparison_1_MM}
	\end{subfigure}
	\hfill
	\begin{subfigure}[b]{0.49\textwidth}
		\centering
		\includegraphics[width=\linewidth]{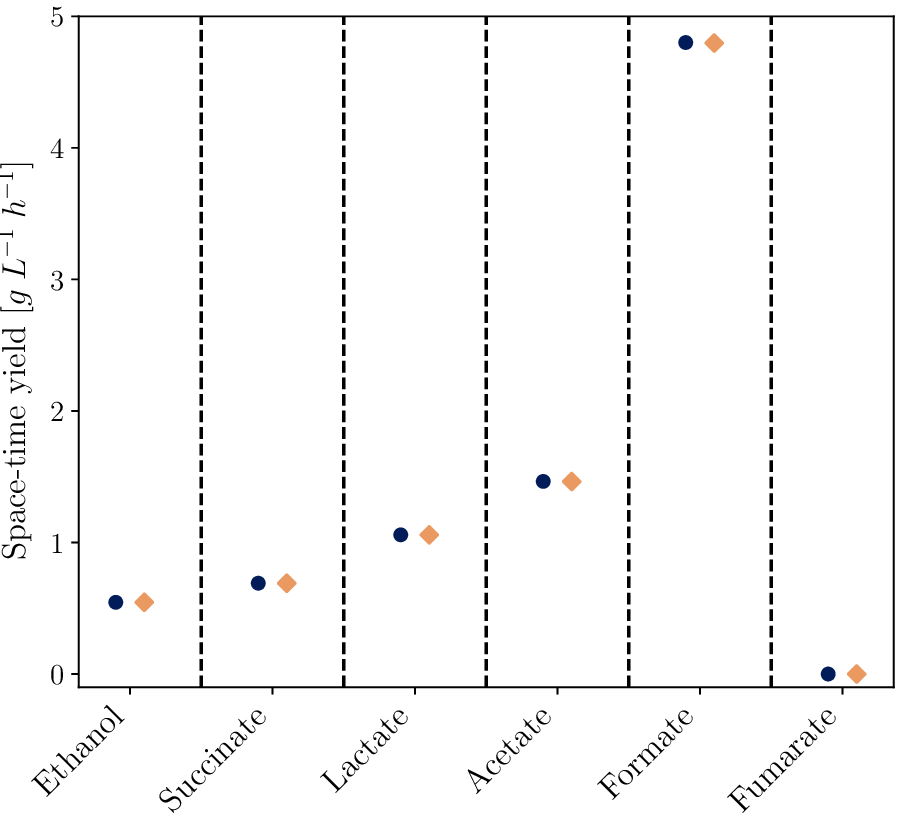}
		\caption{$Monod$}
		\label{subfig:comparison_1_Monod}
	\end{subfigure}
	\hfill
	\caption{Comparison of space-time yield using Michaelis-Menten and Monod kinetics for one reaction elimination on the genome-scale metabolic model \textit{i}ML1515 \citep{Monk.2017}. Ethanol, succinate, and lactate were produced via anaerobic pathways; acetate, formate, and fumarate were produced via aerobic pathways.}
	\label{fig:comparison_1KO}
\end{figure}

\clearpage

\begin{table}[htb]
	\caption{One reaction elimination: Comparison of space-time yield, product concentration, substrate concentration, biomass concentration, growth rate, substrate uptake flux $v_S$, product flux $v_P$, and molar yield using Michaelis-Menten and Monod kinetics for the genome-scale metabolic model \textit{i}ML1515 \citep{Monk.2017}. Note that, for all results, the substrate feed concentration reached its upper bound of \SI{10}{\gram\per\liter}. STY: Space-time  yield, MM: Michaelis-Menten}
	\label{tab:data_comparison_1KO}
	\resizebox{\textwidth}{!}{
		\begin{tabular}{lllllllllllll}
			\hline
			Chemical&Kinetic&Optimization &Knockouts                         &STY &Product &Substrate &Biomass&Growth&v$_{s}$ &v$_{p}$ & Yield \\
			&&&& [\si{\gram\per\liter\per\hour}]&[\si{\gram\per\liter}]&[\si{\gram\per\liter}]&[\si{\gramdryweight\per\liter}]&[\si{\per\hour}]&\multicolumn{2}{c}{[\si{\mmol\per\gramdryweight\per\hour}]}&[\si{\mol\per\mol}]\\\hline
			Ethanol&MM & Sequential &Acetate reversible transport via proton symport (periplasm)&0.43&4.11&1.09&0.56&0.1&9.19&16.57&1.8\\
			& & SimulKnock & Acetate reversible transport via proton symport (periplasm)&0.43&4.12&1.07&0.56&0.1&9.18&16.57&1.8\\
			&Monod & Sequential &Acetate reversible transport via proton symport (periplasm)&0.55&4.58&0.01&0.66&0.12&10.0&17.93&1.79\\
			& & SimulKnock & Acetate reversible transport via proton symport (periplasm)&0.55&4.58&0.01&0.66&0.12&10.0&17.93&1.79\\
			Succinate&MM & Sequential &Ethanol transport via diffusion (extracellular to periplasm)&0.54&5.23&1.1&0.56&0.1&9.2&8.39&0.91\\
			& & SimulKnock & Ethanol transport via diffusion (extracellular to periplasm)&0.54&5.25&1.1&0.56&0.1&9.2&8.42&0.91\\
			&Monod & Sequential &Ethanol transport via diffusion (extracellular to periplasm)&0.69&5.84&0.01&0.66&0.12&10.0&9.07&0.91\\
			& & SimulKnock & Ethanol transport via diffusion (extracellular to periplasm)&0.69&5.84&0.01&0.66&0.12&10.0&9.07&0.91\\
			Lactate&MM & Sequential &Ethanol transport via diffusion (extracellular to periplasm)&0.83&8.01&1.1&0.56&0.1&9.2&16.74&1.82\\
			& & SimulKnock & Ethanol transport via diffusion (extracellular to periplasm)&0.83&8.04&1.09&0.56&0.1&9.19&16.78&1.83\\
			&Monod & Sequential &Ethanol transport via diffusion (extracellular to periplasm)&1.06&8.95&0.01&0.66&0.12&10.0&18.12&1.81\\
			& & SimulKnock & Ethanol transport via diffusion (extracellular to periplasm)&1.06&8.95&0.01&0.66&0.12&10.0&18.12&1.81\\
			Acetate&MM & Sequential &ATP synthase (four protons for one ATP) (periplasm)&1.19&4.68&0.96&1.4&0.25&9.09&14.37&1.58\\
			& & SimulKnock & ATP synthase (four protons for one ATP) (periplasm)&1.19&4.69&0.96&1.4&0.25&9.09&14.38&1.58\\
			&Monod & Sequential &ATP synthase (four protons for one ATP) (periplasm)&1.46&5.14&0.02&1.58&0.28&10.0&15.71&1.57\\
			& & SimulKnock & ATP synthase (four protons for one ATP) (periplasm)&1.46&5.13&0.03&1.58&0.28&10.0&15.71&1.57\\
			Formate&MM & Sequential &ATP synthase (four protons for one ATP) (periplasm)&3.9&15.34&0.96&1.4&0.25&9.1&61.8&6.79\\
			& & SimulKnock & ATP synthase (four protons for one ATP) (periplasm)&3.9&15.35&0.96&1.4&0.25&9.09&61.79&6.79\\
			&Monod & Sequential &ATP synthase (four protons for one ATP) (periplasm)&4.8&16.85&0.02&1.58&0.28&10.0&67.56&6.76\\
			& & SimulKnock & ATP synthase (four protonsfor one ATP) (periplasm)&4.8&16.83&0.03&1.58&0.28&10.0&67.56&6.76\\
			Fumarate&MM & Sequential &Fumarase&0.35&0.48&0.91&4.04&0.72&9.05&0.76&0.08\\
			& & SimulKnock & Fumarase&0.35&0.47&1.15&3.93&0.74&9.24&0.78&0.08\\
			&Monod & Sequential &Tryptophanase (L-tryptophan)&0.0&0.0&0.01&0.0&0.16&10.0&0.0&0.0\\
			& & SimulKnock & 6-phosphogluconolactonase&0.0&0.0&0.01&0.0&0.16&10.0&0.0&0.0\\
			\hline
		\end{tabular}
	}
\end{table}
\clearpage
\subsection{Data for two reaction eliminations}
\begin{table}[htb]
	\caption{Two reaction eliminations: Comparison of space-time yield, product concentration, substrate concentration, biomass concentration, growth rate, substrate uptake flux $v_S$, product flux $v_P$, and molar yield using Michaelis-Menten and Monod kinetics for the genome-scale metabolic model 	\textit{i}ML1515 \citep{Monk.2017}. Note that, for all results, the substrate feed concentration reached its upper bound of \SI{10}{\gram\per\liter}. STY: Space-time yield, MM: Michaelis-Menten}
	\label{tab:data_comparison_2KO}
	\resizebox{\textwidth}{!}{
		\begin{tabular}{llllllllllll}
			\hline
			Chemical&Kinetic&Optimization &Knockouts                         &STY &Product &Substrate &Biomass&Growth&v$_{s}$ &v$_{p}$ & Yield \\
			&&&& [\si{\gram\per\liter\per\hour}]&[\si{\gram\per\liter}]&[\si{\gram\per\liter}]&[\si{\gramdryweight\per\liter}]&[\si{\per\hour}]&\multicolumn{2}{c}{[\si{\mmol\per\gramdryweight\per\hour}]}&[\si{\mol\per\mol}]\\\hline
			Ethanol&MM & Sequential &1. ATP synthase (four protons for one ATP) (periplasm)&0.09&4.2&1.64&0.11&0.02&9.45&18.58&1.97\\
			(anaerob)&&&2. Triose-phosphate isomerase&&&&&&&\\
			& & SimulKnock &1. 6-phosphogluconolactonase&0.43&4.16&1.1&0.56&0.1&9.2&16.81&1.83\\
			&&&2. Acetate reversible transport via proton symport (periplasm)&&&&&&&\\
			&Monod & Sequential &1. ATP synthase (four protons for one ATP) (periplasm)&0.13&5.02&0.0&0.14&0.03&10.0&19.62&1.96\\
			&&&2. Triose-phosphate isomerase&&&&\\
			& & SimulKnock &1. 6-phosphogluconolactonase&0.55&4.65&0.01&0.65&0.12&10.0&18.19&1.82\\
			&&&2. Acetate reversible transport via proton symport (periplasm)&&&&\\
			Succinate&MM & Sequential &1. ATP synthase (four protons for one ATP) (periplasm)&0.42&7.0&1.23&0.32&0.06&9.28&11.51&1.24\\
			(anaerob)&&&2. D-glucose transport via PEP:Pyr PTS (periplasm)&&&&&&&\\
			& & SimulKnock &1. Acetaldehyde dehydrogenase (acetylating)&0.54&5.25&1.09&0.56&0.1&9.2&8.41&0.91\\
			&&&2. D-alanine-D-alanine dipeptidase&&&&&&&\\
			&Monod & Sequential &1. ATP synthase (four protons for one ATP) (periplasm)&0.55&7.98&0.0&0.39&0.07&10.0&12.39&1.24\\
			&&&2. D-glucose transport via PEP:Pyr PTS (periplasm)&&&&\\
			& & SimulKnock &1. Acetaldehyde dehydrogenase (acetylating)&0.69&5.84&0.01&0.66&0.12&10.0&9.07&0.91\\
			&&&2. Cytosine deaminase&&&&\\
			Lactate&MM & Sequential &1. ATP synthase (four protons for one ATP) (periplasm)&0.17&8.12&1.65&0.11&0.02&9.45&18.59&1.97\\
			(anaerob)&&&2. Triose-phosphate isomerase&&&&&&&\\
			& & SimulKnock &1. 6-phosphogluconolactonase&0.83&8.04&1.1&0.56&0.1&9.2&16.8&1.83\\
			&&&2. Acetate reversible transport via proton symport (periplasm)&&&&&&&\\
			&Monod & Sequential &1. ATP synthase (four protons for one ATP) (periplasm)&0.25&9.7&0.0&0.14&0.03&10.0&19.62&1.96\\
			&&&2. Triose-phosphate isomerase&&&&\\
			& & SimulKnock &1. 6-phosphogluconolactonase&1.06&8.98&0.01&0.65&0.12&10.0&18.19&1.82\\
			&&&2. Acetate reversible transport via proton symport (periplasm)&&&&\\
			Acetate&MM & Sequential &1. ATP synthase (four protons for one ATP) (periplasm)&0.55&5.34&1.1&0.55&0.1&9.2&16.85&1.83\\
			(aerob)&&&2. Phosphoglycerate mutase&&&&&&&\\
			& & SimulKnock &1. Phosphoserine transaminase&1.2&4.85&0.96&1.36&0.25&9.1&14.89&1.64\\
			&&&2. ATP synthase (four protons for one ATP) (periplasm)&&&&&&&\\
			&Monod & Sequential &1. ATP synthase (four protons for one ATP) (periplasm)&0.7&5.97&0.01&0.65&0.12&10.0&18.23&1.82\\
			&&&2. Phosphoglycerate mutase&&&&\\
			& & SimulKnock &1. ATP synthase (four protons for one ATP)(periplasm)&1.48&5.32&0.03&1.53&0.28&10.0&16.29&1.63\\
			&&&2. Phosphoglycerate dehydrogenase&&&&\\
			Formate&MM & Sequential &1. ATP synthase (four protons for one ATP) (periplasm)&1.85&17.96&1.1&0.55&0.1&9.2&74.23&8.07\\
			(aerob)&&&2. Phosphoglycerate mutase&&&&&&&\\
			& & SimulKnock &1. Phosphofructokinase (s7p)&3.9&15.35&0.96&1.4&0.25&9.09&61.79&6.79\\
			&&&2. ATP synthase (four protons for one ATP) (periplasm)&&&&&&&\\
			&Monod & Sequential &1. ATP synthase (four protons for one ATP) (periplasm)&2.35&20.06&0.01&0.65&0.12&10.0&80.35&8.04\\
			&&&2. Phosphoglycerate mutase&&&&\\
			& & SimulKnock &1. Phosphogluconate dehydrogenase&4.8&16.83&0.03&1.58&0.28&10.0&67.56&6.76\\
			&&&2. ATP synthase (four protons for one ATP) (periplasm)&&&&\\
			Fumarate&MM & Sequential &1. Fumarase&0.35&0.48&0.9&4.04&0.72&9.04&0.75&0.08\\
			(aerob)&&&2. Cytochrome oxidase bd (menaquinol-8: 2 protons) (periplasm)&&&&&&&\\
			& & SimulKnock &1. Fumarase&0.35&0.49&0.9&4.04&0.72&9.04&0.77&0.08\\
			&&&2. 1-tetradecanoyl-sn-glycerol 3-phosphate O-acyltransferase (n-C14:0)&&&&&&&\\
			&Monod & Sequential &1. Fumarase&0.41&0.51&0.41&4.28&0.8&10.0&0.84&0.08\\
			&&&2. Cytochrome oxidase bd (menaquinol-8: 2 protons) (periplasm)&&&&\\
			& & SimulKnock &1. Fumarase&0.29&0.42&0.78&3.54&0.69&10.0&0.72&0.07\\
			&&&2. Cytochrome oxidase bo3 (ubiquinol-8: 4 protons) (periplasm)&&&&\\
			\hline
		\end{tabular}
	}
\end{table}
\clearpage
\subsection{Data for three reaction eliminations}
\begin{table}[htb]
	\caption{Three reaction eliminations: Comparison of space-time yield, product concentration, substrate concentration, biomass concentration, growth rate, substrate uptake flux $v_S$, product flux $v_P$, and molar yield using Monod kinetics for the genome-scale metabolic model \textit{i}ML1515 \citep{Monk.2017} for aerobic production. Note that, for all results, the substrate feed concentration reached its upper bound of \SI{10}{\gram\per\liter}.  STY: Space-time yield}
	\label{tab:data_comparison_3KO}
	\resizebox{\textwidth}{!}{
		\begin{tabular}{lllllllllll}
			\hline
			Chemical&Optimization &Knockouts                         &STY &Product &Substrate &Biomass&Growth&v$_{s}$ &v$_{p}$ & Yield \\
			&&& [\si{\gram\per\liter\per\hour}]&[\si{\gram\per\liter}]&[\si{\gram\per\liter}]&[\si{\gramdryweight\per\liter}]&[\si{\per\hour}]&\multicolumn{2}{c}{[\si{\mmol\per\gramdryweight\per\hour}]}&[\si{\mol\per\mol}]\\\hline
			Ethanol & Sequential &1. ATP synthase (four protons for one ATP) (periplasm)&0.6&4.6&0.01&0.72&0.13&10.0&17.99&1.8\\
			&&2. Acetate reversible transport via proton symport (periplasm)&&&&&&&\\
			&&3. Succinate dehydrogenase (irreversible)&&&&&&&\\
			& SimulKnock &1. Phosphotransacetylase&0.6&4.45&0.01&0.75&0.14&10.0&17.41&1.74\\
			&&2. ATP synthase (four protons for one ATP) (periplasm)&&&&\\
			&&3. Succinyl-CoA synthetase (ADP-forming)&&&&\\
			Succinate & Sequential &1. Fumarase&0.75&5.86&0.01&0.71&0.13&10.0&9.1&0.91\\
			&&2. ATP synthase (four protons for one ATP) (periplasm)&&&&&&&\\
			&&3. Acetate reversible transport via proton symport (periplasm)&&&&&&&\\
			& SimulKnock &1. Phosphotransacetylase&0.76&5.6&0.01&0.75&0.14&10.0&8.7&0.87\\
			&&2. ATP synthase (four protons for one ATP) (periplasm)&&&&\\
			&&3. Succinyl-CoA synthetase (ADP-forming)&&&&\\
			Lactate & Sequential &1. 2-Oxogluterate dehydrogenase&1.16&8.69&0.01&0.74&0.13&10.0&17.58&1.76\\
			&&2. ATP synthase (four protons for one ATP) (periplasm)&&&&&&&\\
			&&3. Acetate reversible transport via proton symport (periplasm)&&&&&&&\\
			& SimulKnock &1. Phosphotransacetylase&1.16&8.6&0.01&0.75&0.14&10.0&17.41&1.74\\
			&&2. ATP synthase (four protons for one ATP) (periplasm)&&&&\\
			&&3. Succinyl-CoA synthetase (ADP-forming)&&&&\\
			Acetate & Sequential &1. 6-phosphogluconate dehydratase&0.87&7.6&0.01&0.64&0.11&10.0&23.21&2.32\\
			&&2. ATP synthase (four protons for one ATP) (periplasm)&&&&&&&\\
			&&3. Triose-phosphate isomerase&&&&&&&\\
			& SimulKnock &1. Phosphoserine transaminase&2.24&4.1&0.13&3.0&0.55&10.0&12.68&1.27\\
			&&2. Phosphogluconate dehydrogenase&&&&\\
			&&3. Succinate dehydrogenase (irreversible)&&&&\\
			Formate & Sequential &1. 2-dehydro-3-deoxy-phosphogluconate aldolase&2.57&22.32&0.01&0.64&0.11&10.0&89.39&8.94\\
			&&2. ATP synthase (four protons for one ATP) (periplasm)&&&&&&&\\
			&&3. Triose-phosphate isomerase&&&&&&&\\
			& SimulKnock &1. ATP synthase (four protons for one ATP) (periplasm)&4.8&16.83&0.03&1.58&0.28&10.0&67.56&6.76\\
			&&2. Glucosyltransferase II (LPS core synthesis)&&&&\\
			&&3. Glycerol-3-phosphate acyltransferase (C16:0)&&&&\\
			Fumarate & Sequential &1. Fumarase&0.74&5.75&0.01&0.71&0.13&10.0&9.1&0.91\\
			&&2. ATP synthase (four protons for one ATP) (periplasm)&&&&&&&\\
			&&3. Acetate reversible transport via proton symport (periplasm)&&&&&&&\\
			& SimulKnock &1. Fumarase&2.27&4.12&0.14&3.02&0.55&10.0&6.6&0.66\\
			&&2. Phosphoserine transaminase&&&&\\
			&&3. Phosphogluconate dehydrogenase&&&&\\
			\hline
		\end{tabular}
	}
\end{table}

\subsection{Data for Figure \ref{fig:Lab_Comparison}: experimental data compared with SimulKnock results}
\begin{table}[]
	\caption{Comparison of SimulKnock results with lab data from \citet{vanHeerden.2013}.}
	\centering
	\begin{tabular}{llll}
		\hline
		Glucose feed concentration & Dilution rate    &  space-time yield  & source \\
		\si{\gram\per\liter}&\si{\per\hour}&\si{\gram\per\liter\per\hour}& -\\\hline
		20 &0.06&0.72&\\
		&0.15&0.73&\\
		&0.15&0.68&\\
		&0.09&0.85&\citet{vanHeerden.2013}\\
		&0.09&0.87&\\
		&0.19&0.46&\\
		&0.04&0.52&\\\cline{2-4}
		&0.13&1.51&SimulKnock, \textit{i}ML1515\\
		&0.22&2.42&SimulKnock, \textit{i}EC1349\_Crooks\\\hline
		50&0.06&0.71&\\
		&0.06&0.69&\\
		&0.09&0.65&\citet{vanHeerden.2013}\\
		&0.09&0.59&\\
		&0.04&0.58&\\
		&0.04&0.57&\\
		&0.02&0.38&\\\cline{2-4}
		&0.13&3.79&SimulKnock, \textit{i}ML1515\\
		&0.22&6.05&SimulKnock, \textit{i}EC1349\_Crooks\\\hline
	\end{tabular}
	\label{tab:data_plot_experimental}
\end{table}

\end{document}